\documentclass{article}
\usepackage{amsfonts}
\usepackage{amsmath}
\usepackage{graphicx}
\usepackage{subcaption}
\usepackage{float}
\usepackage{booktabs}
\usepackage{longtable}

\newtheorem{theorem}{Theorem}

\newtheorem{conclusion}{Conclusion}

\newtheorem{definition}{Definition}
\newtheorem{example}{Example}

\begin{document}

\title{The multivariate local dependence function}
\author{Ismihan Bayramoglu and Pelin Er\c{s}in \\
%EndAName
Department of Mathematics, Izmir University of Economics}
\maketitle

\begin{abstract}
The local dependence function is important in many applications of probability and statistics. We extend the bivariate local dependence function introduced by Bairamov and Kotz (2000) and further developed by Bairamov et al. (2003) to three-variate and multivariate local dependence function characterizing the dependency between three and more random variables in a given specific point. The definition and properties of the three-variate local dependence function are discussed. An example of a three-variate local dependence function for underlying three-variate normal distribution is presented. The graphs and tables with numerical values are provided. The multivariate extension of the local dependence function that can characterize the dependency between multiple random variables at a specific point is also discussed.

Key Words: Dependent random variables, joint distribution, measures of dependency, local dependence function.
\end{abstract}

\section{Introduction}

The concept of dependence is the most important subject in probability and its applications. In applications of probability and statistics involving two or more random variables, the accurate model can be designed only if
information about the form of dependence is given or assumed. In the case of bivariate random variables, there has been introduced and studied well-known measures of dependence such as Kendall's tau, Spearman's rho measuring the concordance and discordance between random variables, and Pearson's correlation coefficient measuring linear dependence, among others. These measures are given in the form of a scalar, i.e. only one number that describes the
degree of dependence between random variables in the set of all possible values of these random variables. More precisely, assume that $X$ and $Y$ are two random variables taking all values from the set $G\in \mathbb{R}^{2}.$ Then, for instance, Pearson's correlation coefficient $\rho_{X,Y}$
is a number between $-1$ and $1,$ regardless $X$ and $Y$ can take any values from $G,$ i.e. the dependence measured by $\rho _{X,Y}$ is the same for all points of the region $G.$ However, in many applications it may be required
to measure the dependency between random variables $X$ and $Y$ in a given particular point $(x_{0},y_{0})\in G,$ assuming that the dependency at this point may differ from the dependency at another point $(x_{1},y_{1})\in G.$
For example, the dependence between random variables that describe the complex of symptoms in the healthy tissue in the human body may be different than the dependence on the infected tissue. Therefore, there is a need to consider the association measures that can characterize the dependency at
any point of the region $G$. The local dependence measure introduced by Holland and Wang (1987), Bjerve and Doksum (1993)  Doksum et al. (1994), Blyth (1993, 1994 a, b),and Jones (1996) satisfy most of the desirable properties for a measure of dependence described by Scarsini (1984) (see Nelsen (2006)) and also the R\'{e}nyi's axioms (see R\'{e}nyi (1959)). However, these measures have some weaknesses, for example, the local dependence function of Holland and Wang (1987) (Jones (1996))
\begin{equation*}
	H_{1}(x,y)=\frac{\partial ^{2}\log f(x,y)}{\partial x\partial y}
\end{equation*}
for bivariate normal distribution is a constant for all values of $(x,y),$ where, $f(x,y)$ is the joint probability density function of $X$ and $Y.$ Furthermore, the class of all distributions involving an exponential family with its canonical parameter being a linear function has constant local dependence. Jones (1996, 1998) identifies all distributions other than bivariate normal with constant local dependence and mentions that this characterization is not very helpful.

The local dependence function introduced by Bairamov and Kotz (2000) and further developed by Bairamov et al. (2003) is a bivariate function that characterizes the dependence between two random variables at a particular point is based on the localization of the Pearson's correlation coefficient
to the point. The local dependence function possesses the properties of a dependence measure between two random variables. Furthermore, the value of this local dependence function in a particular point equals Pearson's correlation coefficient. Therefore, it can also be regarded as a functional generalization of Pearson's correlation coefficient. Since the local dependence function is a bivariate function, it can be used to form the local dependence maps for bivariate random variables, which describe the degree of dependence between random variables in the given area. Bairamov-Kotz local dependence function has other interesting and attractive properties beyond satisfying the R\'{e}ny's axioms and overcomes weaknesses
of previous measures of local dependence measures. Bairamov-Kotz local dependence function for random variables $X$ and $Y$ is defined as
\begin{equation*}
	H(x,y)=\frac{E(X-E(X\mid Y=y))(Y-E(Y\mid X=x))}{\sqrt{E(X-E(X\mid Y=y))^{2}}%
		\sqrt{E(Y-E(Y\mid X=x))^{2}}},
\end{equation*}
and can also be shown as
\begin{equation*}
	H(x,y)=\frac{\rho _{X,Y}+\phi _{X}(y)\phi _{Y}(x)}{\sqrt{1+\phi _{Y}^{2}(x)}%
		\sqrt{1+\phi _{X}^{2}(y)}},
\end{equation*}
where
\begin{eqnarray*}
	\phi _{X}(y) &=&\frac{EX-E(X\mid Y=y)}{\sigma _{X}} \\
	\phi _{Y}(x) &=&\frac{EY-E(Y\mid X=x)}{\sigma _{Y}} \\
	\sigma _{X}^{2} &=&Var(X),\sigma _{Y}^{2}=Var(Y).
\end{eqnarray*}
The method of constructing $H(x,y)$ was based on the idea of replacing $EX$ with its best predictor $E(X\mid Y=y)$ of $X$ through $Y$ at the point $y,$ and replacing $Y$ with its the best predictor $E(Y\mid X=x)$ of $Y$ through $X,$ at the point $x$ in the expression of Pearson's correlation coefficient.
This replacement would detect the degree of dependence of $X$ and $Y$ at the point $(x,y).$ The most intriguing property of $H(x,y)$ is that if, in its expression, we replace $t=\phi _{X}(y),$ $y=\phi _{X}^{-1}(t)$ and $
s=\phi _{Y}(x),$ $x=\phi _{Y}^{-1}(s),$ where $\phi_{X}^{-1}(t)$ and $\phi_{Y}^{-1}(s)$ are the generalized inverse functions of $\phi _{X}(y)$ and $\phi _{Y}(x),$ respectively, then
\begin{equation*}
	h(t,s)=\frac{\rho _{X,Y}+ts}{\sqrt{1+t^{2}}\sqrt{1+s^{2}}}
\end{equation*}
has a saddle point at $t=0$ and $s=0,$ with $h(0,0)=\rho_{X,Y},$ which corresponds to points $x^{\ast }$ and $y^{\ast }$ such that $\phi_{Y}(x^{\ast })=0$ and $\phi _{X}(y^{\ast })=0$ $,$ i.e. $EX=E(X\mid Y=y^{\ast })$ and $EY=E(Y\mid X=x^{\ast }).$ It is clear that,
\begin{equation*}
	H(x,y)=h(\phi _{X}(y),\phi _{Y}(x))
\end{equation*}
and
\begin{equation*}
	h(t,s)=H(\phi _{Y}^{-1}(s),\phi _{X}^{-1}(t)).
\end{equation*}
Note that the above considerations are true if the functions $\phi _{X}(y)$ and $\phi _{Y}(x)$ are monotone.

Unlike Holland-Wang-Jones local dependence function $H_{1}(x,y)$, Bairamov-Kotz local dependence function for bivariate normal distribution with zero mean and unit standard deviation is not constant. It is
\begin{eqnarray*}
H(x,y) &=&\frac{\rho _{X,Y}+xy\rho_{X,Y}^{2}}{\sqrt{1+x^{2}\rho _{X,Y}^{2}}\sqrt{1+y^{2}\rho _{X,Y}^{2}}} \\
&=&h(x\rho _{X,Y},y\rho _{X,Y}),
\end{eqnarray*}
indicating that $H(x,y)$ is a more powerful measure of dependence than $H_{1}(x,y)$ (for more discussions and graphical comparisons see Kotz and Nadarajah (2003)).

In this work, we extend the Bairamov-Kotz local dependence function to three-variate and multivariate random variables. The definition and properties of the three-variate local dependence function are discussed.
Some examples of particular distributions are presented. The paper is organized as follows: Section 2 provides preliminary information and discussion on the bivariate local dependence function. In Section 3, the new
three-variate local dependence function based on the conditional expectation concept is introduced and the properties are discussed. An example of a three-variate local dependence function for underlying three-variate normal distribution is given. The graphical representations and numerical values calculated for selected parameters and particular points are given in tables. In Section 4, we discuss the multivariate extension of the local dependence function to any finite number of parameters.

\section{Preliminaries. The Bivariate Local Dependence Function}

The local dependence measure introduced by Holland and Wang (1987), Bjerve and Doksum (1993), Doksum et al. (1994), Blyth (1993, 1994 a,b),and Jones (1996) satisfy most of the desirable properties for a measure of dependence described by Scarsini (1984) (see Nelsen (2006)) and also the R\'{e}nyi's axioms (see R\'{e}nyi (1959)).

However, these measures have some weaknesses, for example the local dependence function of Holland and Wang (1987) and (Jones (1996))
\begin{equation*}
	H_{1}(x,y)=\frac{\partial ^{2}\log f(x,y)}{\partial x\partial y}
\end{equation*}
for bivariate normal distribution is a constant for all values of $(x,y),$ where, $f(x,y)$ is the joint probability density function of $X$ and $Y.$ Furthermore, the class of all distributions involving an exponential family with its canonical parameter being a linear function has constant local dependence. Jones (1996, 1998) identifies all distributions other than bivariate normal with constant local dependence and mentions that this characterization is not very helpful.

The local dependence function introduced by Bairamov and Kotz (2000) and further developed by Bairamov et al. (2003) is a bivariate function that characterizes the dependence between two random variables at a particular point. The local dependence function possesses the properties of a dependence measure between two random variables.

Furthermore, the value of this local dependence function in a particular point equals Pearson's correlation coefficient. Therefore, it can also be regarded as a functional generalization of Pearson's correlation coefficient.

Since the local dependence function is a bivariate function, it can be used to form the local dependence maps for bivariate random variables, which describes the degree of dependence between random variables in the given area. For example, the local dependence maps in medicine are visualization tools that provide insights into how individual patient features (such as demographics, lab results, and medical history) influence model predictions
for specific medical outcomes. A nice application of local dependence maps constructed with the Bairamov-Kotz local dependence function using real clinical data is given in Oru\c{c} and Hudaverdi (2009).

A measure of local dependence introduced in Bairamov and Kotz (2000) (see also Bairamov, Kotz, and Kozubowski (2003)) is based on the localization of Pearson's correlation coefficient to the point. Bairamov-Kotz local dependence function (see Bairamov and Kotz (2000)) has other interesting and attractive properties beyond satisfying the R\'{e}ny's axioms and overcomes weaknesses of previous measures of local dependence measures.

Bairamov-Kotz local dependence function for random variables $X$ and $Y$ is defined as
\begin{equation*}
	H(x,y)=\frac{E(X-E(X\mid Y=y))(Y-E(Y\mid X=x))}{\sqrt{E(X-E(X\mid Y=y))^{2}}%
		\sqrt{E(Y-E(Y\mid X=x))^{2}}},
\end{equation*}

and can also be shown as
\begin{equation*}
	H(x,y)=\frac{\rho _{X,Y}+\phi _{X}(y)\phi _{Y}(x)}{\sqrt{1+\phi _{Y}^{2}(x)}%
		\sqrt{1+\phi _{X}^{2}(y)}},
\end{equation*}%
where
\begin{eqnarray*}
	\phi _{X}(y) &=&\frac{EX-E(X\mid Y=y)}{\sigma _{X}} \\
	\phi _{Y}(x) &=&\frac{EY-E(Y\mid X=x)}{\sigma _{Y}} \\
	\sigma _{X}^{2} &=&Var(X),\sigma _{Y}^{2}=Var(Y).
\end{eqnarray*}

The method of constructing $H(x,y)$ was based on the idea of replacing $EX$ with its best predictor $E(X\mid Y=y)$ of $X$ through $Y$ at the point $y,$ and replacing $Y$ with its the best predictor $E(Y\mid X=x)$ of $Y$ through $X,$at the point $x$ in the expression of Pearson's correlation coefficient. This replacement would detect the degree of dependence of $X$ and $Y$ at the point $(x,y).$

The most intriguing property of $H(x,y)$ is that if in its expression we replace $t=\phi _{X}(y),$ $y=\phi_{X}^{-1}(t)$ and $s=\phi _{Y}(x),$ $x=\phi _{Y}^{-1}(s),$ where $\phi _{X}^{-1}(t)$ and $\phi _{Y}^{-1}(s)$ are the generalized inverse functions of $\phi _{X}(y)$ and $\phi _{Y}(x),$ respectively, then
\begin{equation*}
	h(t,s)=\frac{\rho _{X,Y}+ts}{\sqrt{1+t^{2}}\sqrt{1+s^{2}}}
\end{equation*}%
has a saddle point at $t=0$ and $s=0,$ with $h(0,0)=\rho _{X,Y},$ which corresponds to points $x^{\ast }$ and $y^{\ast }$ such that $\phi_{Y}(x^{\ast })=0$ and $\phi _{X}(y^{\ast })=0$ $,$ i.e. $EX=E(X\mid Y=y^{\ast })$ and $EY=E(Y\mid X=x^{\ast }).$ 
It is clear that,
\begin{equation*}
	H(x,y)=h(\phi _{X}(y),\phi _{Y}(x))
\end{equation*}
and%
\begin{equation*}
	h(t,s)=H(\phi _{Y}^{-1}(s),\phi _{X}^{-1}(t)).
\end{equation*}
Note that the above considerations are true if the functions $\phi _{X}(y)$ and $\phi _{Y}(x)$ are monotone.

Unlike Holland-Wang-Jones local dependence function $H_{1}(x,y)$, Bairamov-Kotz local dependence function for bivariate normal distribution with zero mean and unit standard deviation is not constant, it is
\begin{eqnarray*}
	H(x,y) &=&\frac{\rho _{X,Y}+xy\rho _{X,Y}^{2}}{\sqrt{1+x^{2}\rho _{X,Y}^{2}}%
		\sqrt{1+y^{2}\rho _{X,Y}^{2}}} \\
	&=&h(x\rho _{X,Y},y\rho _{X,Y}),
\end{eqnarray*}%
indicating that $H(x,y)$ is more powerfull measure of dependence than $H_{1}(x,y)$(for more discussions and graphical comparisons see Kotz and Nadarajah (2003)).

Kotz and Nadarajah (2003) provide new motivations for local dependence function introduced by Bairamov and Kotz (2000) which is a localized version of Galton-Pearson correlation coeefficient. 
\bigskip

\subsection{Properties of local dependence function H(x,y)}

The properties of $H(x,y)$ are given in the following theorem:

\begin{theorem}
\label{Theorem 1}(Bairamov and Kotz (2000)) \textit{The local dependence function }$H(x,y)$\textit{\ has the following properties:}

\textit{1}$^{\circ }$\textit{. If }$X$\textit{\ and }$Y$\textit{\ are
	independent then }$H(x,y)=0$\textit{\ for any }$(x,y)\in N_{X,Y}.$

\textit{2}$^{\circ }$\textit{. }$\left\vert H(x,y)\right\vert \leq 1,$%
\textit{\ for alll }$(x,y)\in N_{X,Y}.$

\textit{3}$^{\circ }$\textit{. If }$\left\vert H(x,y)\right\vert =1$\textit{%
	\ for some }$(x,y)\in N_{X,Y}$\textit{\ then }$\rho \neq 0.$

\textit{4}$^{\circ }.$\textit{\ Let }$H(x,y)=0$\textit{\ for all }$(x,y)\in
N_{X,Y}.$\textit{\ Then one of the following assertions are valid.}

$a)$\textit{\ Either }$EX$\textit{\ }$=E(X\mid Y=y)$\textit{\ , }$E(Y\mid
X=x)=const$\textit{\ or }$EY$\textit{\ }$=E(Y\mid X=x)$\textit{\ , }$E(X\mid
Y=y)=const$\textit{\ for all }$x$\textit{\ and }$y$\textit{\ and }$\rho =0.$

$b)$\textit{\ Either }$EX=E(X\mid Y=y)$\textit{\ for all }$y,$\textit{\ }$%
E(Y\mid X=x)$\textit{\ is not a constant function or }$EY=E(Y\mid X=x)$%
\textit{\ for all }$x,$\textit{\ }$E(X\mid Y=y)$\textit{\ is not a constant
	function and }$\rho =0.$

$c)$\textit{\ }$E(X\mid Y=y)=const,$\textit{\ }$E(Y\mid X=x)=const$\textit{\
	for all }$x$\textit{\ and }$y,$\textit{\ i.e. }$X$\textit{\ and }$Y$\textit{%
	\ are regression independent.}

\textit{5}$^{\circ }.$\textit{\ Let }$\left\vert \rho \right\vert =1$\textit{%
	\ and assume that }$H(x,y)=1$\textit{\ at the point }$(x,y)$\textit{\ then }$%
\varphi _{X}(y)=\varphi _{Y}(x)$\textit{\ up to a sign.}

\textit{6}$^{\circ }.$\textit{\ }$h(t,z)$\textit{\ has a saddle point at }$%
(0,0).$\textit{\ That is the point }$(x^{\ast },y^{\ast })$\textit{\
	satisfying }$\varphi _{X}(y^{\ast })=\varphi _{Y}(x^{\ast })=0$\textit{\ is
	a saddle point of }$H(x,y)$\textit{\ and }$H$\textit{\ }$(x^{\ast },y^{\ast
})=\rho .$\textit{\ }

\end{theorem}

\subsection{Examples}

\textbf{1.} Consider a bivariate normal distribution with joint p.d.f. 
\begin{equation*}
f(x,y)=\frac{1}{2\pi \sigma _{X}\sigma _{Y}\sqrt{1-\rho^{2}}}\exp \left\{ -\frac{1}{2(1-\rho ^{2})}\left( \frac{x^{2}}{\sigma _{X}^{2}}-2\rho \frac{xy}{\sigma _{X}\sigma _{Y}}+\frac{y^{2}}{\sigma _{Y}^{2}}\right) \right\} ,
\end{equation*}
where $\sigma _{X}^{2}=Var(X),$ $\sigma _{Y}^{2}=Var(Y),$ $EX=0,EY=0,$ $\rho=\frac{Cov(X,Y)}{\sigma _{X}\sigma _{Y}}.$
\\

$E(Y\mid X=x)=\rho \frac{\sigma _{Y}}{\sigma _{X}}x,$ $E(X\mid Y=y)=\rho \frac{\sigma _{X}}{\sigma _{Y}}y$ and 
\begin{equation*}
H(x,y)=\frac{\rho +\rho ^{2}xy}{\sqrt{\sigma _{X}^{2}+\rho^{2}\frac{\sigma_{X}^{2}}{\sigma _{Y}^{2}}y^{2}}\sqrt{\sigma _{Y}^{2}+\rho ^{2}\frac{\sigma_{Y}^{2}}{\sigma _{X}^{2}}x^{2}}}.
\end{equation*}
In the simple case when $\sigma _{X}=\sigma _{Y}=1$ we have 
\begin{equation}
H(x,y)=\frac{\rho +\rho ^{2}xy}{\sqrt{1+\rho ^{2}y^{2}}\sqrt{1+\rho ^{2}x^{2}}} \label{bn1}
\end{equation}

\thinspace \thinspace \thinspace

\begin{figure}[H]
    \centering
    \begin{subfigure}[t]{0.45\textwidth}
        \centering
        \includegraphics[width=\textwidth]{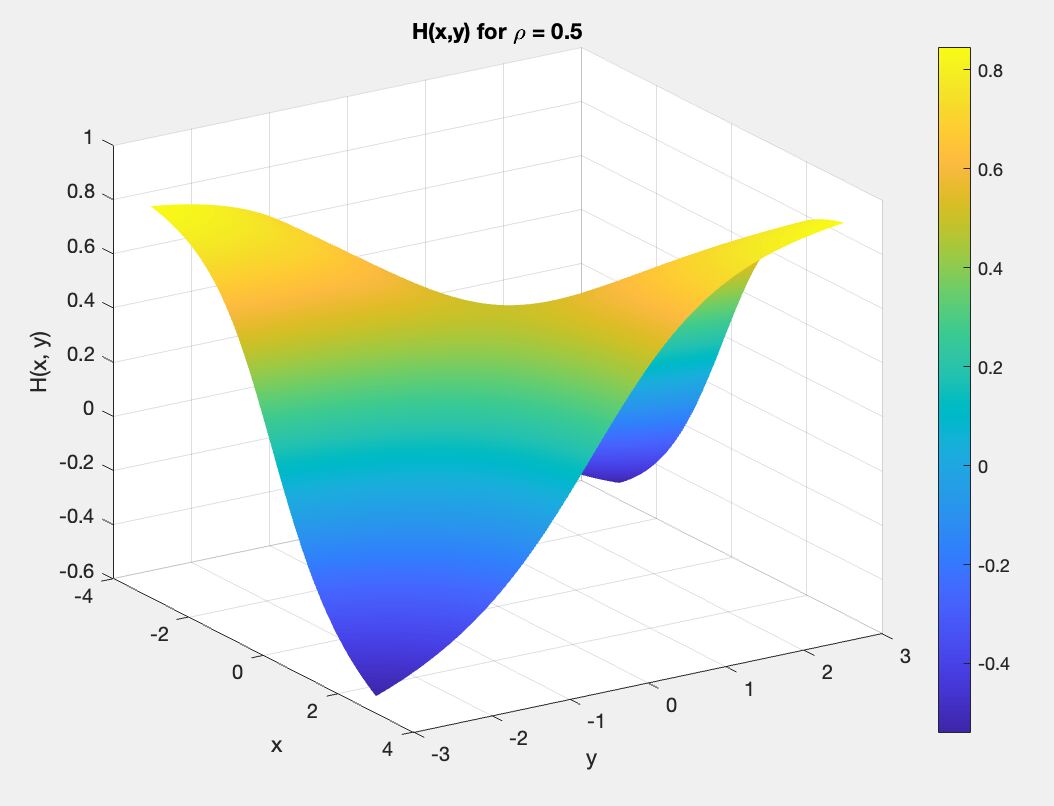}
    \end{subfigure}
    \hfill
    \begin{subfigure}[t]{0.45\textwidth}
        \centering
        \includegraphics[width=\textwidth]{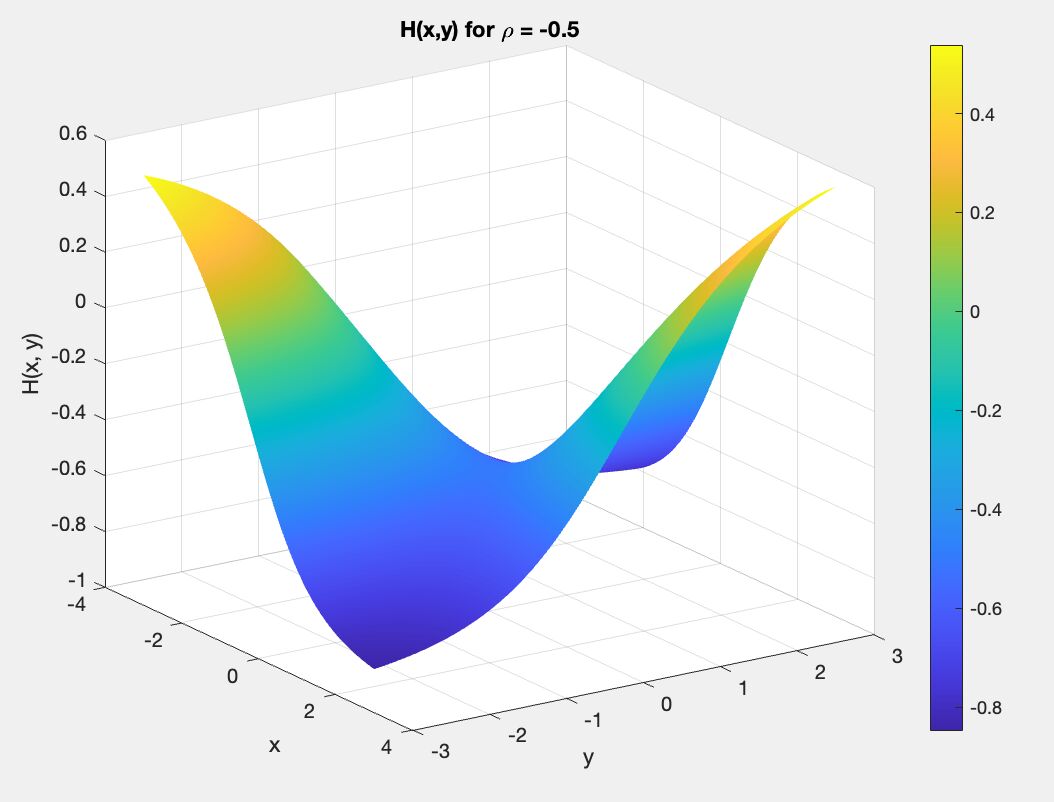}
    \end{subfigure}
    \caption{Graph of \( H(x,y) \) given in (\ref{bn1}) for \( \rho = 0.5 \) and \( \rho = -0.5 \), with \( -3 \leq x, y \leq 3 \).}
    \label{fig:1}
\end{figure}

\begin{figure}[H]
    \centering
    \begin{subfigure}[t]{0.45\textwidth}
        \centering
        \includegraphics[width=\textwidth]{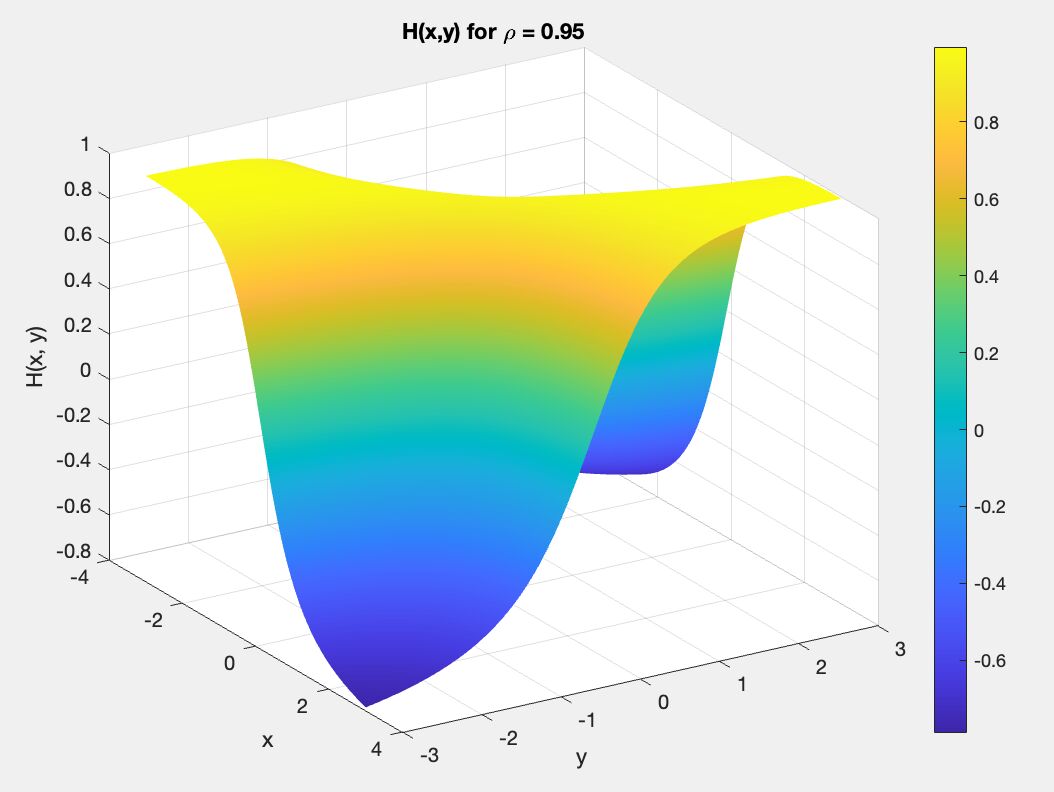}
    \end{subfigure}
    \hfill
    \begin{subfigure}[t]{0.45\textwidth}
        \centering
        \includegraphics[width=\textwidth]{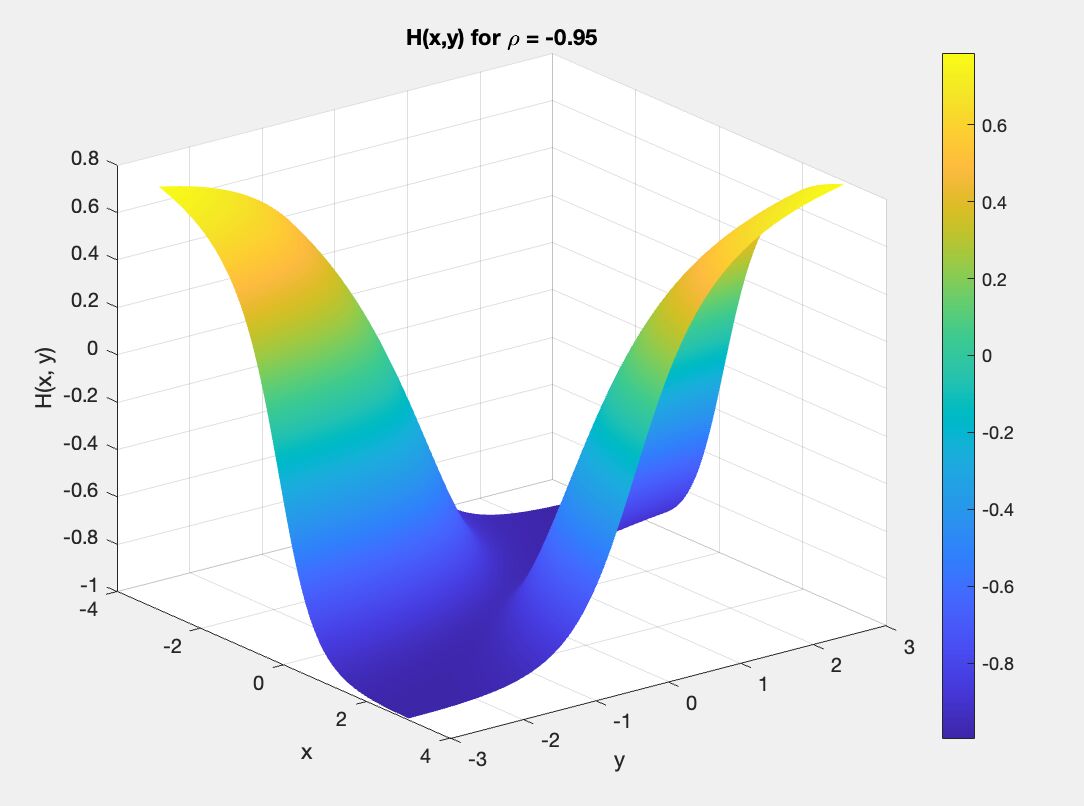}
    \end{subfigure}
    \caption{Graph of \( H(x,y) \) given in (\ref{bn1}) for \( \rho = 0.95 \) and \( \rho = -0.95 \), with \( -3 \leq x, y \leq 3 \).}
    \label{fig:2}
\end{figure}

\section{The three-variate extension of local dependence function}

In many practical applications, a local dependence function is needed to characterize local dependence among multiple random variables. Applying the same idea used in the construction of the Bairamov-Kotz local dependence function $H(x,y),$ we aim to introduce the three-variate and then the multivariate reasonable local dependence function which satisfies desirable
properties similar to R\'{e}ny's axioms in bivariate case, and characterizes the dependence among these random variables in a given particular point. For these reasons, first, we consider the random variables $X,Y,Z$ with joint distribution function 
\begin{equation*}
F_{X,Y,X}(x,y,z)=C(F_{X}(x),F_{Y}(y),F_{Z}(z)),
\end{equation*}
where $C$ is a connecting copula and $F_{X}(x),F_{Y}(y)$ and $F_{Z}(z)$ are the marginal distribution functions of $X,Y$ and $Z,$ respectively. Let the Pearson's correlation coefficients between the random variables $(X,Y),$ $(X,Z)$ and $(Y,Z)$ be 
\begin{equation*}
\rho _{X,Y}=\frac{Cov(X,Y)}{\sigma _{X}\sigma _{Y}},\rho _{X,Z}=\frac{Cov(X,Z)}{\sigma _{X}\sigma _{Z}},\rho _{Y,Z}=\frac{Cov(Y,Z)}{\sigma_{Y}\sigma _{Z}},
\end{equation*}
where 
\begin{eqnarray*}
Cov(X,Y) &=&E(X-EX)(Y-EY) \\
Cov(X,Z) &=&E(X-EX)(Z-EZ) \\
Cov(Y,Z) &=&E(Y-EY)(Z-EZ)
\end{eqnarray*}
and $\sigma _{X}^{2}=Var(X),\sigma _{Y}^{2}=Var(Y)$ and $\sigma_{Z}^{2}=Var(Z).$

\bigskip

Consider the function
\begin{equation*}
\xi_{X}(y,z)=EX-E(X\mid Y=y,Z=z)\newline
,
\end{equation*}
where $E(X\mid Y=y,Z=z)$ is the conditional expectation of $X,$ given $Y=y,Z=z.$ Similarly, denote 
\begin{align*}
\xi_{Y}(x,z)& =EY-E(Y\mid X=x,Z=z)\newline
, \\
\xi_{Z}(x,y)& =EZ-E(Z\mid X=x,Y=y).
\end{align*}
It is well known that the conditional expectation of $X,$ given $Y$ and $Z,$ defined as $E(X\mid Y,Z)\equiv L(Y,Z),$ where $L(y,z)=E(X\mid Y=y,Z=z),$ is the best predictor for $X$ through the random variables $Y$ and $Z,$ in the sense of the mean squares distance. More precisely, 
\begin{equation*}
\underset{a(y,z)\in A(y,z)}{\min }E(X-a(Y,Z))^{2}=E(X-L(Y,Z))^{2},
\end{equation*}
where $A(y,z)$ is the class of all bivariate measurable functions and the minimum is taken over all bivariate measurable functions $a(y,z),(y,z)\in \mathbb{R}^{2}.$ Therefore, replacing $EX$ with $E(X\mid Y=y,Z=z),$ we expect to
estimate the value of $EX$ with given local values of $Y$ and $Z.$ Denote also, 
\begin{align}
\varphi_{X}(y,z)& =\frac{\xi_{X}(y,z)}{\sigma_{X}}=\frac{EX-E(X\mid Y=y,Z=z)\newline}{\sigma_{X}},\newline \notag \\
\varphi_{Y}(x,z)& =\frac{\xi_{Y}(x,z)}{\sigma_{Y}}=\frac{EY-E(Y\mid X=x,Z=z)\newline}{\sigma _{Y}}  \notag \\
\varphi_{Z}(x,y)& =\frac{\xi_{Z}(x,y)}{\sigma_{Z}}=\frac{EZ-E(Z\mid X=x,Y=y)}{\sigma _{Z}}  \label{0}
\end{align}

Consider the functions 
\begin{eqnarray}
&&B(x,y,z)  \notag \\
&\equiv &E[(X-E(X\mid Y=y,Z=z)  \notag \\
\times (Y-E(Y &\mid &X=x,Z=z)  \notag \\
\times (Z-E(Z &\mid &X=x,Y=y)]  \label{01}
\end{eqnarray}
and 
\begin{equation*}
Q(x,y,z)
\end{equation*}
\begin{eqnarray}
&\equiv &E[X-E(X\mid Y=y,Z=z)]^{2}  \notag \\
\times E[Y-E(Y &\mid &X=x,Z=z)]^{2}  \notag \\
\times E[Z-E(Z &\mid &X=x,Y=y)]^{2}.  \label{02}
\end{eqnarray}

Now consider the three-variate function $H(x,y,z),(x,y,z)\in \mathbb{R}^{3}$ defined as follows:
\begin{equation}
H(x,y,z)=\frac{B(x,y,z)}{\sqrt{Q(x,y,z)}}.  \label{1}
\end{equation}
It is easy to see that, 
\begin{eqnarray*}
&&B(x,y,z) \\
&=&E[(X-E(X\mid Y=y,Z=z)(Y-E(Y\mid X=x,Z=z)(Z-E(Z\mid X=x,Y=y)] \\
&=&E[(X-EX+EX-E(X\mid Y=y,Z=z)
\end{eqnarray*}
\begin{equation*}
\times (Y-EY+EY-E(Y\mid X=x,Z=z)
\end{equation*}
\begin{equation*}
\times (Z-EZ+EZ-E(Z\mid X=x,Y=y)]
\end{equation*}

\begin{align}
&= E\left[(X - E[X] + \xi_{X}(y, z)) \right. \notag \\
&\quad \times \left( Y - E[Y] + \xi_{Y}(x, z) \right) \notag \\
&\quad \times \left. \left( Z - E[Z] + \xi_{Z}(x, y) \right) \right] \notag \\
&= E\left[(X - E[X])(Y - E[Y])(Z - E[Z])\right] \notag \\
&\quad + \operatorname{Cov}(Y, Z) \xi_{X}(y, z) + \operatorname{Cov}(X, Z) \xi_{Y}(x, z) \notag \\
&\quad + \operatorname{Cov}(X, Y) \xi_{Z}(x, y) + \xi_{X}(y, z) \xi_{Y}(x, z) \xi_{Z}(x, y).
\label{2}
\end{align}

We also have
\begin{align}
E\left( X - E(X \mid Y = y, Z = z) \right)^{2} &= E[X^{2}] - 2E[X] E(X \mid Y = y, Z = z) \notag \\
&\quad + \left( E(X \mid Y = y, Z = z) \right)^{2} \notag \\
&= E[X^{2}] - (E[X])^{2} + (E[X])^{2} - 2E[X] E(X \mid Y = y, Z = z) \notag \\
&\quad + \left( E(X \mid Y = y, Z = z) \right)^{2} \notag \\
&= \operatorname{Var}(X) + (E[X] - E(X \mid Y = y, Z = z))^{2} \notag \\
&= \operatorname{Var}(X) + \xi_{X}^{2}(y, z) = \sigma_{X}^{2} + \xi_{X}^{2}(y, z).
\label{3}
\end{align}

where 
\begin{equation*}
\xi_{X}^{2}(y,z)=(EX-E(X\mid Y=y,Z=z))^{2}
\end{equation*}
Similarly, 
\begin{align*}
E(Y-E(Y& \mid X=x,Z=z))^{2} \\
& =\sigma_{Y}^{2}+\xi_{Y}^{2}(x,z)
\end{align*}
where 
\begin{eqnarray*}
E(Z-E(Y &\mid &X=x,Y=y))^{2}\newline\\
&=&\sigma _{Z}^{2}+\xi _{Z}^{2}(x,y).
\end{eqnarray*}
Therefore, taking into account (\ref{0}), (\ref{2}) and (\ref{3}) and dividing numerator and denominator of (\ref{1}) by $\sigma_{X},\sigma_{Y},\sigma_{Z}$ we have from (\ref{1})
\begin{eqnarray}
&&H(x,y,z)  \notag \\
&=&\frac{B(x,y,z)}{\sqrt{Q(x,y,z)}}  \notag \\
&=&\frac{\rho _{X,Y,Z}+\rho _{X,Y}\varphi _{Z}(x,y)+\rho_{Y,Z}\varphi_{X}(y,z)+\rho_{X,Z}\varphi_{Y}(x,z)+\varphi_{X}(y,z)\varphi_{Y}(x,z)\varphi_{Z}(x,y)}{\sqrt{1+\varphi_{X}^{2}(y,z)}\sqrt{1+\varphi_{Y}^{2}(x,z)}\sqrt{1+\varphi_{Z}^{2}(x,y)}},  \label{1.1}
\end{eqnarray}
where 
\begin{equation}
\rho _{X,Y,Z}=\frac{E(X-EX)(Y-EY)(Z-EZ)}{\sigma _{x}\sigma _{y}\sigma _{z}}.
\label{1.01}
\end{equation}
Using now (\ref{1.1}) and taking into account (\ref{1.01}) we can give the following definition.

\begin{definition}
The function $H(x,y,z),(x,y,z)\in \mathbb{R}^{3}$ defined as 
\begin{eqnarray*}
&&H(x,y,z) \\
&=&\frac{\rho_{X,Y,Z}+\rho_{X,Y}\varphi_{Z}(x,y)+\rho_{Y,Z}\varphi_{X}(y,z)+\rho_{X,Z}\varphi _{Y}(x,z)+\varphi_{X}(y,z)\varphi_{Y}(x,z)\varphi _{Z}(x,y)}{\sqrt{1+\varphi_{X}^{2}(y,z)}\sqrt{1+\varphi_{Y}^{2}(x,z)}\sqrt{1+\varphi_{Z}^{2}(x,y)}}
\end{eqnarray*}
is called the three-variate total local dependence function for random variables $X,Y$ and $\ Z.$
\end{definition}

Consider the function 
\begin{equation}
h(t,s,w)=\frac{\rho_{X,Y,Z}+t\rho_{X,Y}+s\rho_{Y,Z}+w\rho_{X,Z}+tsw}{\sqrt{1+t^{2}}\sqrt{1+s^{2}}\sqrt{1+w^{2}}}.  \label{1.2}
\end{equation}
It is clear that the function $h(t,s,w)$ is a symmetric function with respect to its variables $t,s$ and $w$. Then replacing $t=\varphi_{Z}(x,y),s=\varphi_{X}(y,z)$ and $w=\varphi_{Y}(x,z)$, we can write 
\begin{equation}
H(x,y,z)=h(\varphi_{Z}(x,y),\varphi_{X}(y,z),\varphi_{Y}(x,z)).
\label{1.3}
\end{equation}
Assume that the system of equations 
\begin{equation*}
\left\{ 
\begin{array}{c}
\begin{array}{c}
t=\varphi_{Z}(x,y) \\ 
s=\varphi_{X}(y,z)
\end{array}
\\ 
w=\varphi_{Y}(x,z)
\end{array}
\right.
\end{equation*}
has a solution with respect to $x,y$ and $z,$ i.e. 
\begin{equation*}
\left\{ 
\begin{array}{c}
\begin{array}{c}
x=\psi_{1}(t,s,w) \\ 
y=\psi_{2}(t,s,w)
\end{array}\\ 
y=\psi_{3}(t,s,w)
\end{array}
\right. ,
\end{equation*}
for some functions $\psi_{1},\psi_{2}$ and $\psi_{3}.$ Then
\begin{equation}
h(t,s,w)=H(\psi_{1}(t,s,w),\psi_{2}(t,s,w),\psi_{3}(t,s,w)).  \label{1.4}
\end{equation}

\bigskip

\begin{theorem}
(Properties of local dependence function) The local dependence function $
H(x,y,z)$ has the following properties: \newline
1. If $X,Y$ and $Z$ are jointly independent random variables, i.e. $C(t,s,w)=tsw,$ $0\leq t,s,w\leq 1,$ then $H(x,y,z)=0$ for all $(x,y,z)\in \mathbb{R}^{3}.$ \newline
2. $\left\vert H(x,y,z)\right\vert \leq 1,$ for all $(x,y,z)\in \mathbb{R}^{3}$ \newline
3. If the random variables $X,Y,Z$ are pairwise independent, i.e. pairs $X$ and $Y$ are independent, $X$ and $Z$ are independent, and $Y$ and $Z$ are independent, then 
\begin{equation*}
H(x,y,z)=\frac{\rho_{X,Y,Z}+\varphi_{X}(y,z)\varphi_{Y}(x,z)\varphi_{Z}(x,y)}{\sqrt{1+\varphi_{X}^{2}(y,z)}\sqrt{1+\varphi_{Y}^{2}(x,z)}\sqrt{1+\varphi_{Z}^{2}(x,y)}}
\end{equation*}
for all $(x,y,z)\in \mathbb{R}^{3}.$\newline
4. \ Let $(x_{1},y_{1},z_{1})\in \mathbb{R}^{3}$ satisfy the system of equations

\begin{equation}
\left\{ 
\begin{array}{c}
EX=E(X\mid Y=y_{1},Z=z_{1}) \\ 
EY=E(Y\mid X=x_{1},Z=z_{1}) \\ 
EZ=E(Z\mid X=x_{1},Y=y_{1}).
\end{array}
\right.  \label{c1}
\end{equation}
\ then $H(x_{1},y_{1},z_{1})=\rho _{X,Y,Z}.$\newline
\end{theorem}

\begin{example}
Consider an $n-$variate random vector with joint normal distribution with p.d.f 
\begin{equation*}
f(\mathbf{x}) = \frac{1}{(2\pi)^{\frac{n}{2}} |\sum|^{\frac{1}{2}}} \exp \left( - (\mathbf{x} - \mathbf{\mu})^{T} \sum^{-1} (\mathbf{x} - \mathbf{\mu}) \right)
\end{equation*}
where
\begin{eqnarray*}
\mathbf{x} &\mathbf{=}&(x_{1},x_{2},....,x_{n}), \\
\mathbf{\mu } &=&(\mu _{1},\mu _{2},...,\mu _{n})^{T} \\
\sum &=&[\sigma_{ij}],i,j=1,2,...,n \\
\left\vert \sum \right\vert &=&\det \left( \sum \right) \\
\mu _{i} &=&EX_{i},\sigma_{ij}=Cov(X_{i},X_{j}).
\end{eqnarray*}
For $n=3$ we have the three-variate random vector $(X,Y,Z)$ with with p.d.f
\begin{equation*}
f(\mathbf{x}\mid \boldsymbol{\mu },\Sigma )=\frac{1}{(2\pi )^{3/2}|\Sigma|^{1/2}}\exp \left( -\frac{1}{2}(\mathbf{x}-\boldsymbol{\mu })^{T}\Sigma^{-1}(\mathbf{x}-\boldsymbol{\mu })\right)
\end{equation*}
and in open form 
\begin{eqnarray*}
f(x_{1},x_{2},x_{3} &\mid &\boldsymbol{\mu },\Sigma ) \\
&=&\frac{1}{(2\pi )^{\frac{3}{2}}\left\vert \sum \right\vert ^{\frac{1}{2}}}\left( \exp -\frac{1}{2\left\vert \sum \right\vert }\left[ a_{11}(x_{1}-\mu_{1})^{2}\right. \right. \\
&&+a_{22}(x_{2}-\mu _{2})^{2}+a_{33}(x_{3}-\mu_{3})^{2}+2a_{12}(x_{1}-\mu_{1})(x_{2}-\mu _{2}) \\
&&+2a_{13}(x_{1}-\mu _{1})(x_{3}-\mu _{3})+\left. \left. +2a_{23}(x_{2}-\mu_{2})(x_{3}-\mu _{3})\right] _{\substack{ \\ }}\right) ,
\end{eqnarray*}
where{}
\begin{equation*}
\left\vert \sum \right\vert =\sigma _{11}\sigma _{22}\sigma_{33}+2\sigma_{12}\sigma _{23}\sigma _{31}-\sigma_{13}^{2}\sigma _{22}-\sigma_{12}^{2}\sigma _{33}-\sigma _{23}^{2}\sigma _{11}
\end{equation*}
$\sum =[\sigma _{ij}]$ is a covariance matrix 
\begin{equation*}
\Sigma =
\begin{bmatrix}
\sigma _{X}^{2} & cov(X,Y) & Cov(X,Z) \\ 
Cov(Y,X) & \sigma _{Y}^{2} & Cov(Y,Z) \\ 
Cov(Z,X) & Cov(Z,Y) & \sigma _{z}^{2}%
\end{bmatrix}
\end{equation*}
and elements of $\sum^{-1}=[a_{ij}]$ are 
\begin{eqnarray*}
a_{11} &=&\sigma_{22}\sigma_{33}-\sigma_{23}^{2} \\
a_{22} &=&\sigma_{11}\sigma_{33-}-\sigma_{13}^{2} \\
a_{33} &=&\sigma_{11}\sigma_{22}-\sigma_{12}^{2} \\
a_{12} &=&\sigma_{13}\sigma_{23}-\sigma_{12}\sigma_{33} \\
a_{13} &=&\sigma_{12}\sigma_{23}-\sigma_{13}\sigma_{22} \\
a_{23} &=&\sigma_{13}\sigma_{12}-\sigma_{11}\sigma_{23}
\end{eqnarray*}
\end{example}

We have visualized the total local dependency $H(x,y,z)$ using a three-variate normal distribution with a fixed covariance matrix and mean vector. We examined how the total local dependency changes when one variable is fixed at different values while allowing the other variables to vary over specified intervals.

The parameters for\textbf{\ }the mean vector and $\mathbf{\Sigma }$ the covariance matrix were initialized as the first step to compute the total local dependency. Then, a grid for two variables was created while the other fixed at a constant value. The following formula was used for the calculation of conditional means of $X,Y,Z$ given the other variables: 
\begin{equation*}
E(X|Y=y,Z=z)=\mu_{X}+\Sigma_{x,(y,z)}\Sigma_{(y,z),(y,z)}^{-1}%
\begin{bmatrix}
y-\mu_{Y} \\ 
z-\mu_{Z}
\end{bmatrix}
\end{equation*}
where, $\Sigma_{X,(Y,Z)}$ is the vector of covariances between $X$ and $(Y,Z)$, 
\begin{equation*}
\Sigma _{X,(Y,Z)}=
\begin{bmatrix}
Cov(X,Y), & Cov(X,Z)
\end{bmatrix}
\end{equation*}

$\Sigma_{(Y,Z),(Y,Y)}$ is the covariance matrix for $Y$ and $Z$, 
\begin{equation*}
\Sigma_{(Y,Z),(Y,Z)}=
\begin{bmatrix}
\sigma_{Y}^{2} & Cov(Y,Z) \\ 
Cov(Z,Y) & \sigma_{Z}^{2}
\end{bmatrix}
\end{equation*}

$\Sigma_{(Y,Z),(Y,Z)}^{-1}$ is the inverse of the covariance matrix of $Y$ and $Z$ which is defined above.

Similar logic was applied for $E(Y|X=x,Z=z)$ and $E(Z|X=x,Y=y)$. Later on, the residuals which represent deviations from the conditional means were computed as they are needed in the total local dependence calculation: 
\begin{equation*}
\xi_{X}(y,z)=\mu_{X}-E(X|Y=y,Z=z)
\end{equation*}
\begin{equation*}
\xi_{Y}(x,z)=\mu_{Y}-E(Y|X=x,Z=z)
\end{equation*}
\begin{equation*}
\xi_{Z}(x,y)=\mu_{Z}-E(Z|X=x,Y=y)
\end{equation*}

Then, to obtain the standardized residuals and normalized covariance standard deviations $\sigma_{x}$, $\sigma_{y}$, and $\sigma_{z}$ as well as the covariance between the variables such as $\text{Cov}(x,y)$ were calculated and determined.

As the final step, the total local dependence formula was applied as follows with a loop that calculates $H(x,y,z)$ for each combination of the changing variables in the grid while keeping the other variable fixed at some value: 
\begin{eqnarray*}
&&H(x,y,z) \\
&=&\frac{\rho_{X,Y,Z}+\rho_{X,Y}\varphi_{Z}(x,y)+\rho_{Y,Z}\varphi_{X}(y,z)+\rho_{X,Z}\varphi_{Y}(x,z)+\varphi _{X}(y,z)\varphi_{Y}(x,z)\varphi_{Z}(x,y)}{\sqrt{1+\varphi_{X}^{2}(y,z)}\sqrt{1+\varphi_{Y}^{2}(x,z)}\sqrt{1+\varphi_{Z}^{2}(x,y)}}
\end{eqnarray*}

where

\begin{equation*}
\rho_{X,Y,Z}=\frac{E(X-EX)(Y-EY)(Z-EZ)}{\sigma_{X}\sigma_{Y}\sigma_{Z}}
\end{equation*}

\begin{eqnarray*}
\varphi_{X}(y,z) &=&\frac{\xi_{X}(y,z)}{\sigma_{X}} \\
\varphi_{Y}(x,z) &=&\frac{\xi_{Y}(x,z)}{\sigma_{Y}} \\
\varphi_{Z}(x,y) &=&\frac{\xi_{Z}(x,y)}{\sigma_{Z}}
\end{eqnarray*}
For the multivariate normal distribution, parameters such as the mean vector and the covariance matrix were considered as follows:

Mean Vector ($\mu $): The mean vector used for the variables $X,Y,Z$: 
\begin{equation}
\mathbf{\mu }=
\begin{bmatrix}
\mu_{X} \\ 
\mu_{Y} \\ 
\mu_{Z}
\end{bmatrix}
=
\begin{bmatrix}
0 \\ 
0 \\ 
0
\end{bmatrix}
\label{nn1}
\end{equation}
We use this mean vector as it centers the distribution at the origin so that the variables are not biased in any particular direction.

Covariance Matrix ($\Sigma $): 
\begin{equation}
\Sigma =
\begin{bmatrix}
\sigma_{X}^{2} & Cov(X,Y) & Cov(X,Z) \\ 
Cov(Y,X) & \sigma_{y}^{2} & Cov(Y,Z) \\ 
Cov(Z,X) & Cov(Z,Y) & \sigma_{z}^{2}
\end{bmatrix}
=
\begin{bmatrix}
1.0 & 0.5 & 0.3 \\ 
0.5 & 1.0 & 0.4 \\ 
0.3 & 0.4 & 1.0
\end{bmatrix}
\label{nn2}
\end{equation}
This matrix defines the variance for each variable ($\sigma_{x}^{2}=\sigma_{y}^{2}=\sigma_{z}^{2}=1$) and the covariance between pairs of variables. The covariance values are taken arbitrarily.

To visualize the total local dependence in a 3D setting, we fixed $x $, $y $, and $z $ at specific values in separate cases and varied the others over specific intervals.

Fixed values for $x $ were $-2$, $0$, and $2$ to analyze how the total local dependency $H(x, y, z) $ changes at these points. Interval for $y $ and $z $ was varied over $y, z \in [-4, 4]$. Additionally, in a separate computation, the interval for $y $ and $z $ was extended to $y, z \in [-100,100] $ to observe the behavior of $H(x, y, z) $ over a larger range.

The following figures show the 3D representation of the total local dependency $H(x,y,z)$ for different fixed values of $x$ and varying intervals for $y$ and $z$.
\begin{figure}[H]
    \centering
    \begin{subfigure}[t]{0.45\textwidth}
        \centering
        \includegraphics[width=\textwidth]{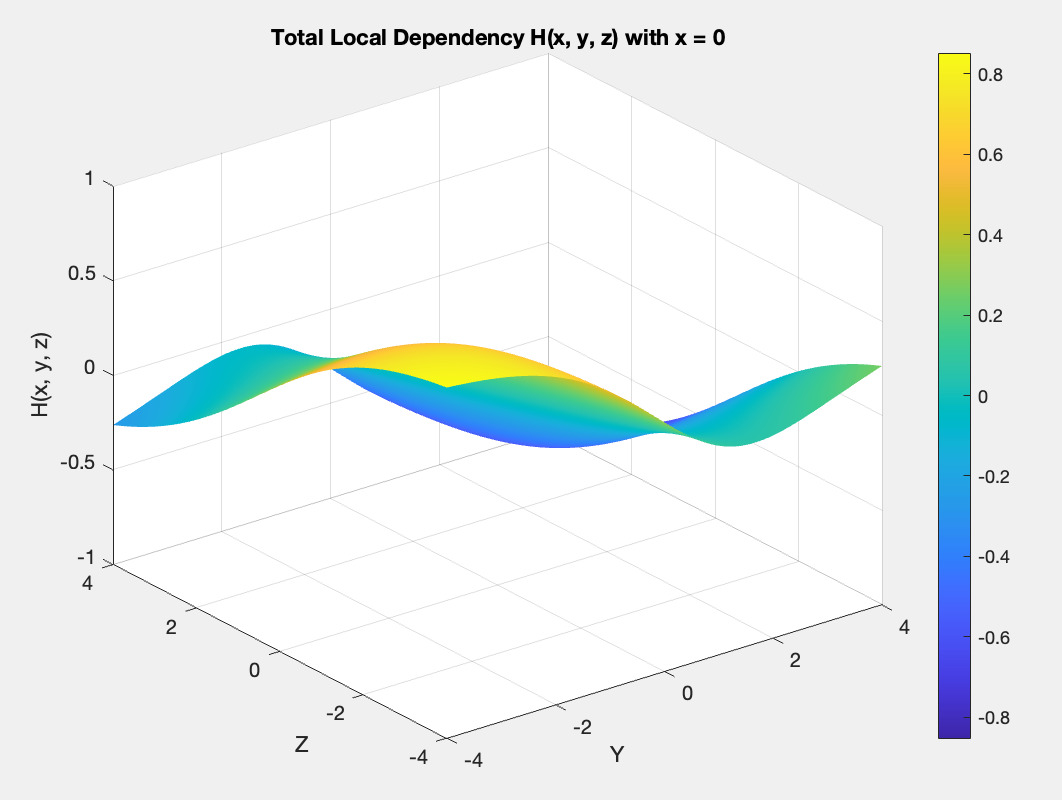}
    \end{subfigure}
    \hfill
    \begin{subfigure}[t]{0.45\textwidth}
        \centering
        \includegraphics[width=\textwidth]{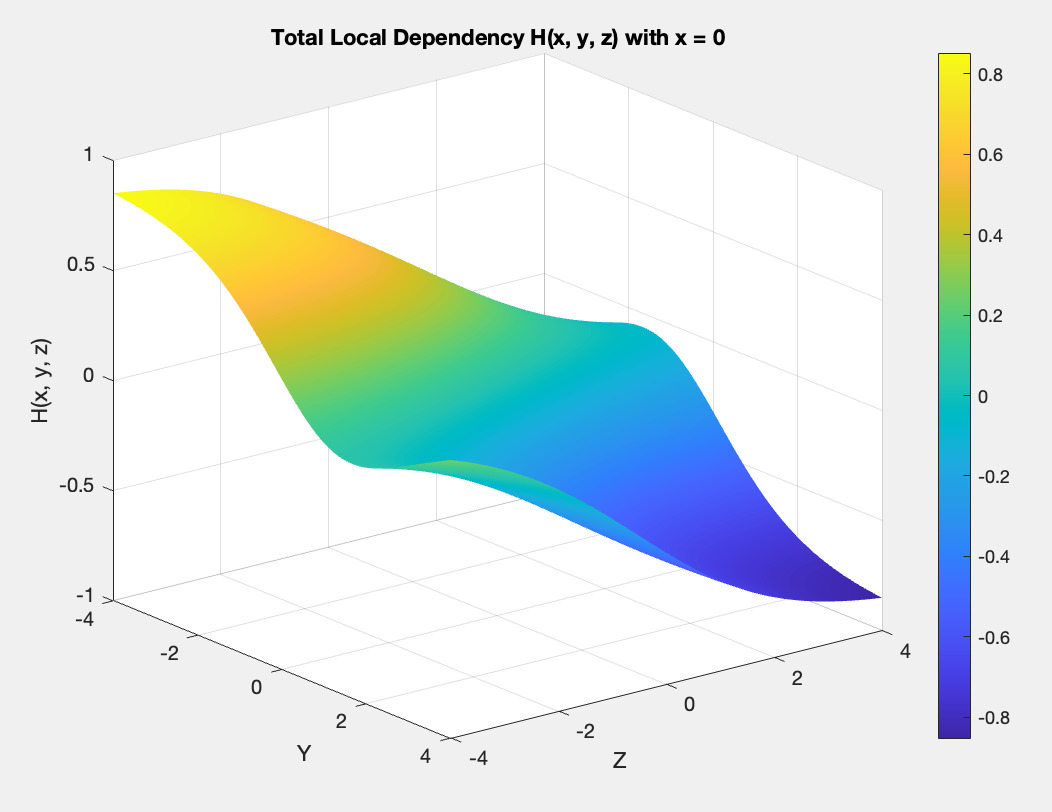}
    \end{subfigure}
    \caption{Graph of \( H(x,y,z) \) for \( x = 0 \),  \( y,z\in (-4,4)\).}
    \label{fig:3}
\end{figure}

\begin{figure}[H]
    \centering
    \begin{subfigure}[t]{0.45\textwidth}
        \centering
        \includegraphics[width=\textwidth]{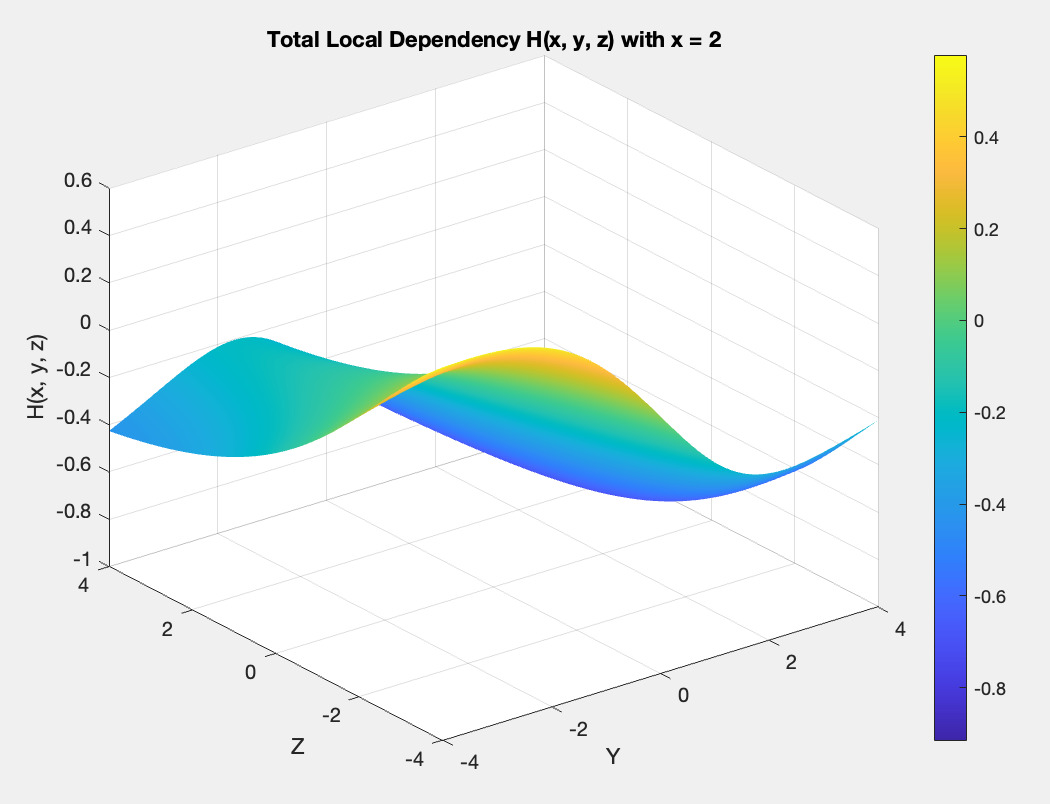}
    \end{subfigure}
    \hfill
    \begin{subfigure}[t]{0.45\textwidth}
        \centering
        \includegraphics[width=\textwidth]{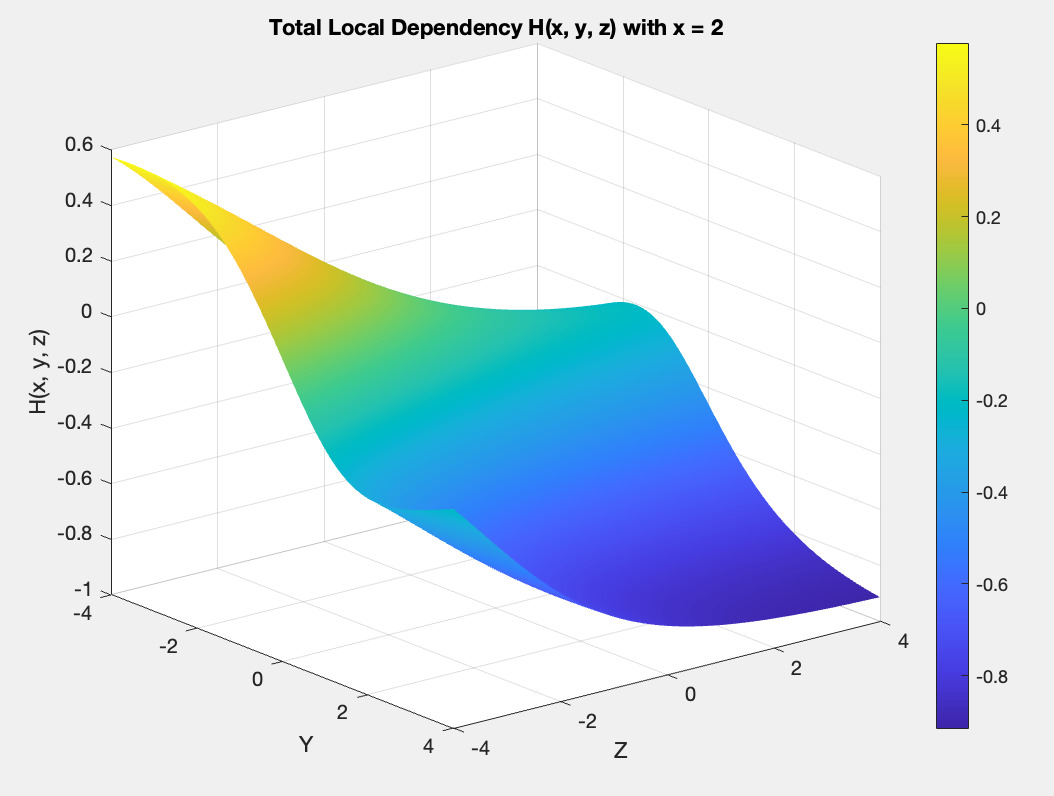}
    \end{subfigure}
    \caption{Graph of \( H(x,y,z) \) for \( x = 2 \),  \( y,z\in (-4,4)\).}
    \label{fig:4}
\end{figure}

\begin{figure}[H]
    \centering
    \begin{subfigure}[t]{0.45\textwidth}
        \centering
        \includegraphics[width=\textwidth]{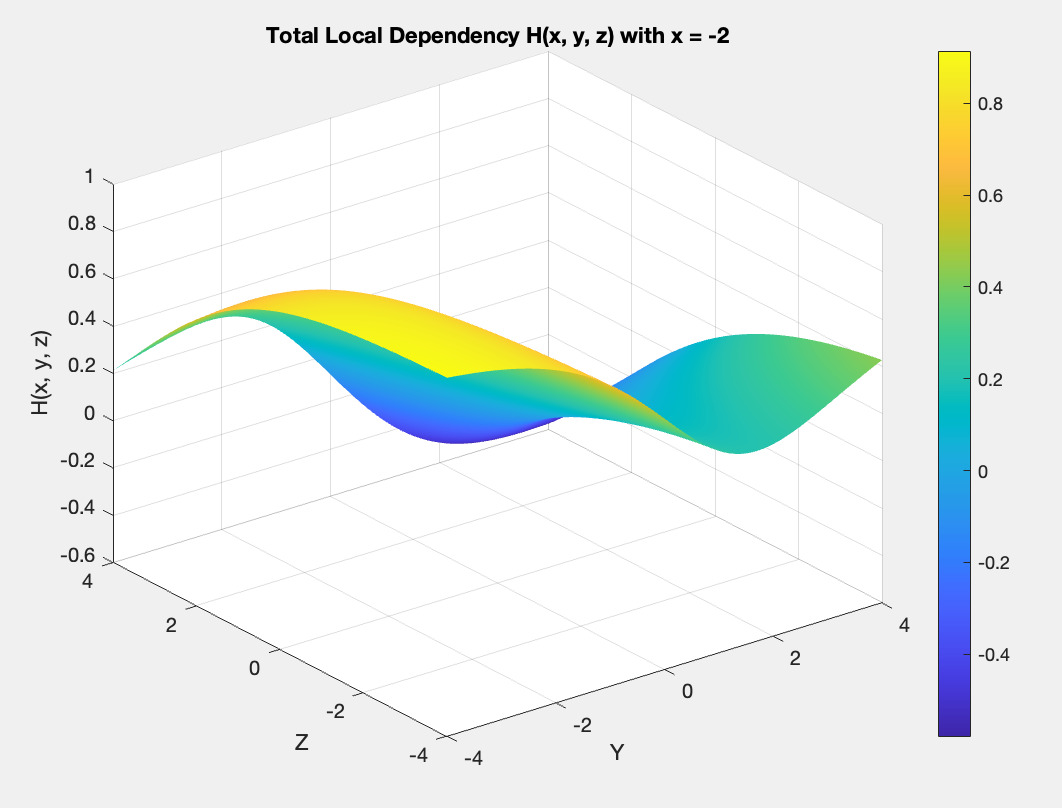}
    \end{subfigure}
    \hfill
    \begin{subfigure}[t]{0.45\textwidth}
        \centering
        \includegraphics[width=\textwidth]{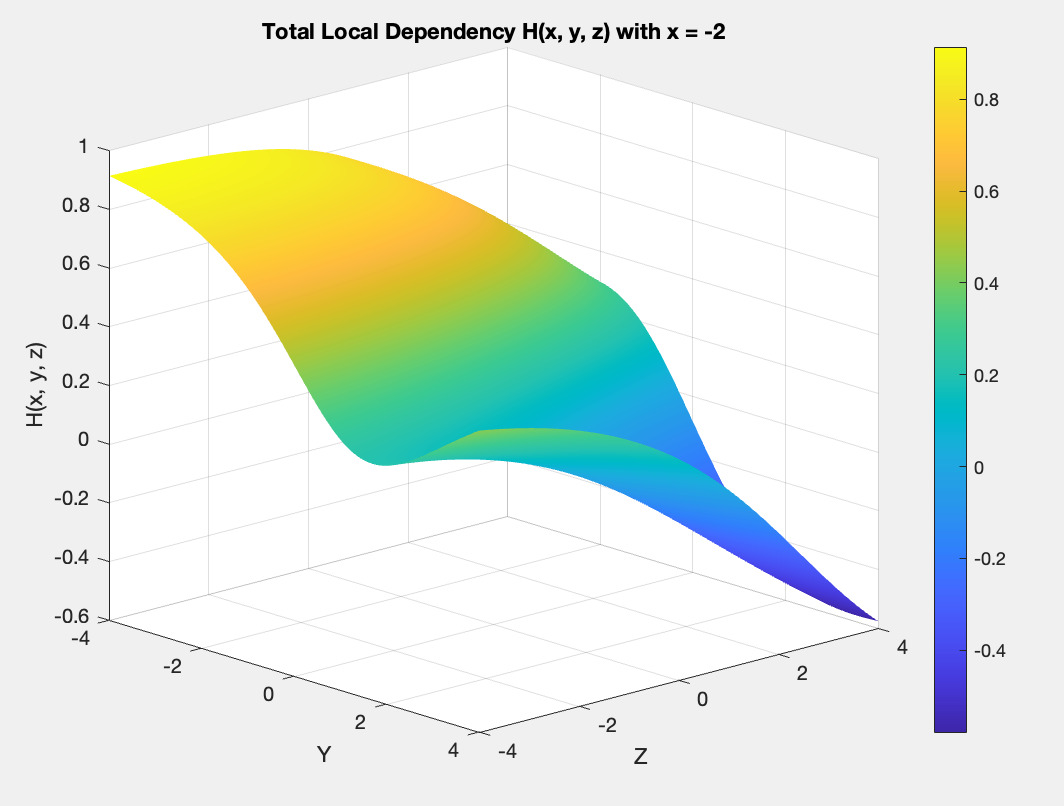}
    \end{subfigure}
    \caption{Graph of \( H(x,y,z) \) for \( x = -2 \),  \( y,z\in (-4,4)\).}
    \label{fig:5}
\end{figure}

\begin{figure}[H]
    \centering
    \begin{subfigure}[t]{0.45\textwidth}
        \centering
        \includegraphics[width=\textwidth]{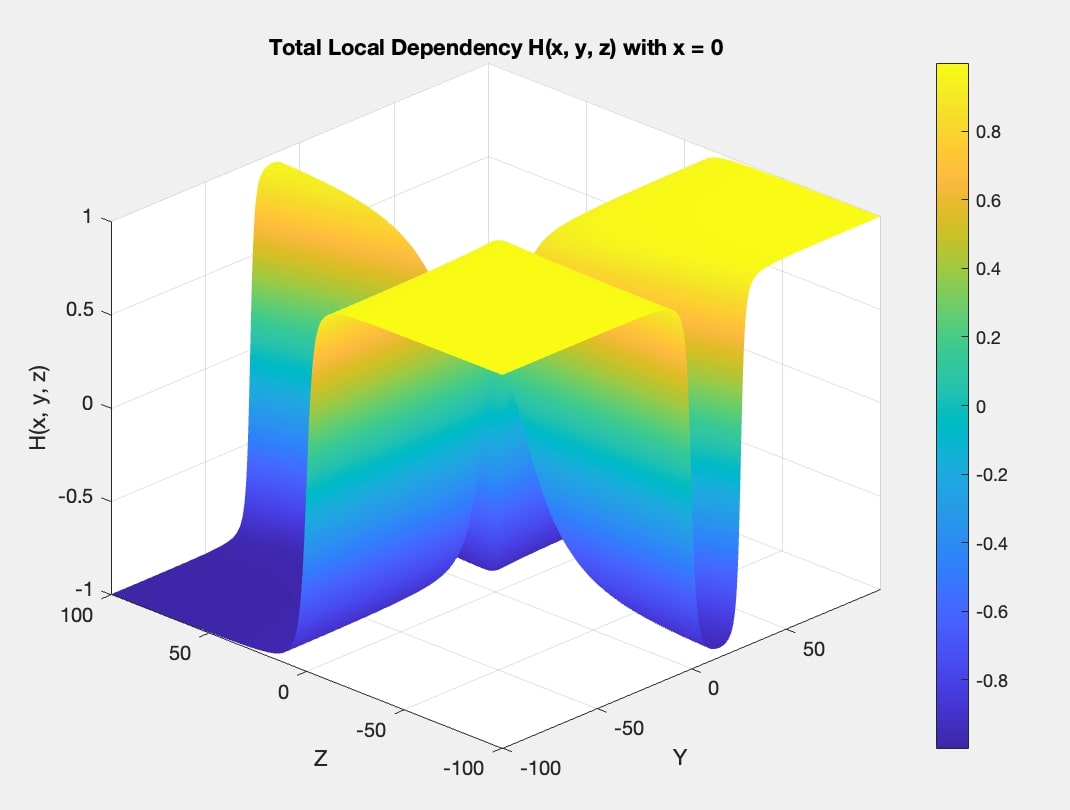}
    \end{subfigure}
    \hfill
    \begin{subfigure}[t]{0.45\textwidth}
        \centering
        \includegraphics[width=\textwidth]{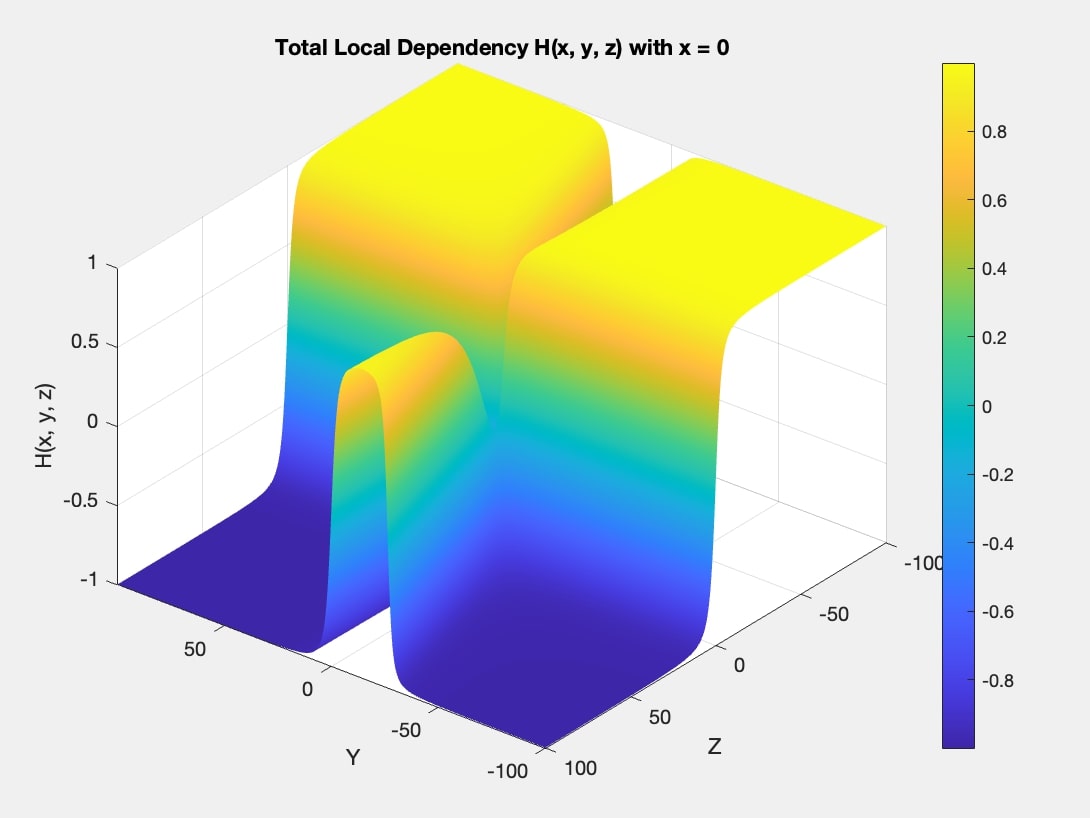}
    \end{subfigure}
    \caption{Graph of \( H(x,y,z) \) for \( x = 0 \),  \( y,z\in (-100,100)\).}
    \label{fig:6}
\end{figure}

In a similar approach to the latter, the 3D representations of total local dependence $H(x,y,z)$ for different fixed values of $y$ and varying intervals for $x$ and $z$ are given.

\begin{figure}[H]
    \centering
    \begin{subfigure}[t]{0.45\textwidth}
        \centering
        \includegraphics[width=\textwidth]{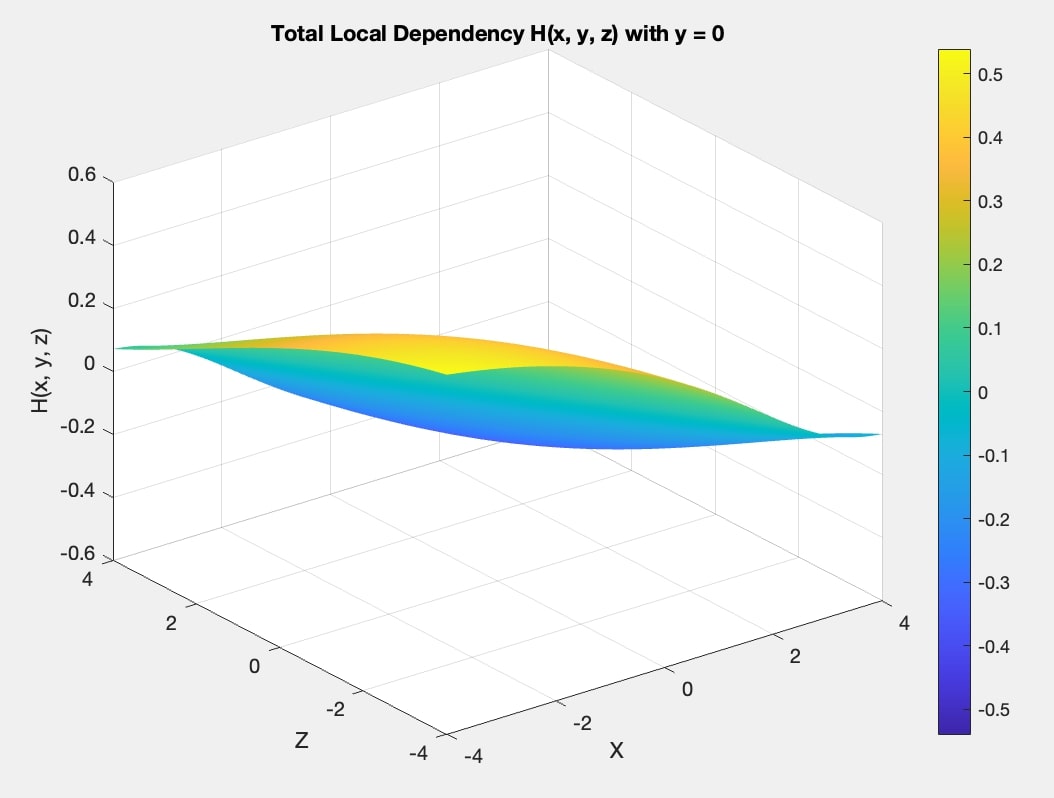}
    \end{subfigure}
    \hfill
    \begin{subfigure}[t]{0.45\textwidth}
        \centering
        \includegraphics[width=\textwidth]{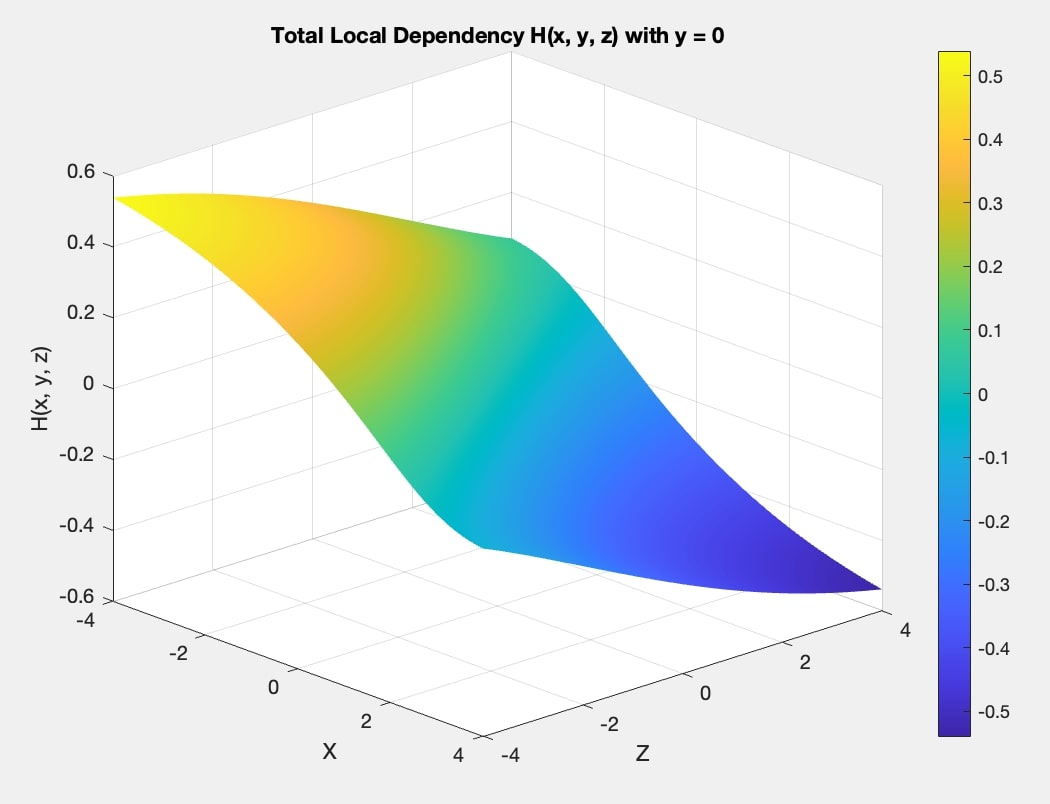}
    \end{subfigure}
    \caption{Graph of \( H(x,y,z) \) for \( y = 0 \),  \(x,z\in (-4,4)\).}
    \label{fig:7}
\end{figure}

\begin{figure}[H]
    \centering
    \begin{subfigure}[t]{0.45\textwidth}
        \centering
        \includegraphics[width=\textwidth]{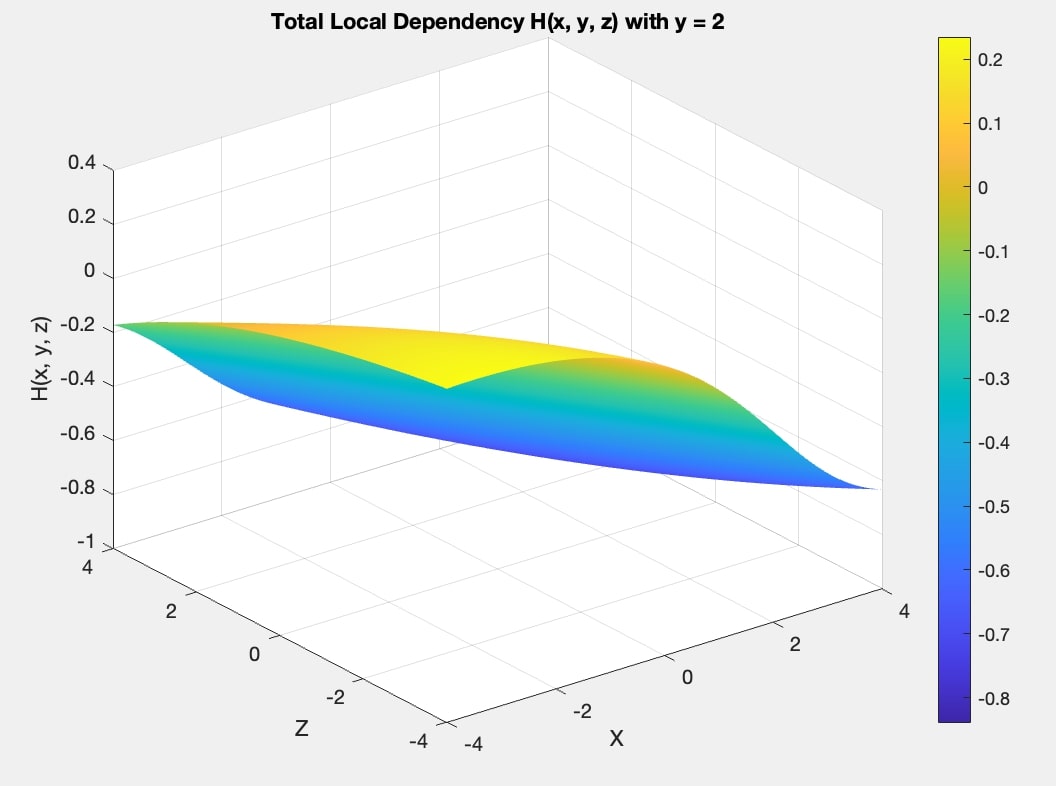}
    \end{subfigure}
    \hfill
    \begin{subfigure}[t]{0.45\textwidth}
        \centering
        \includegraphics[width=\textwidth]{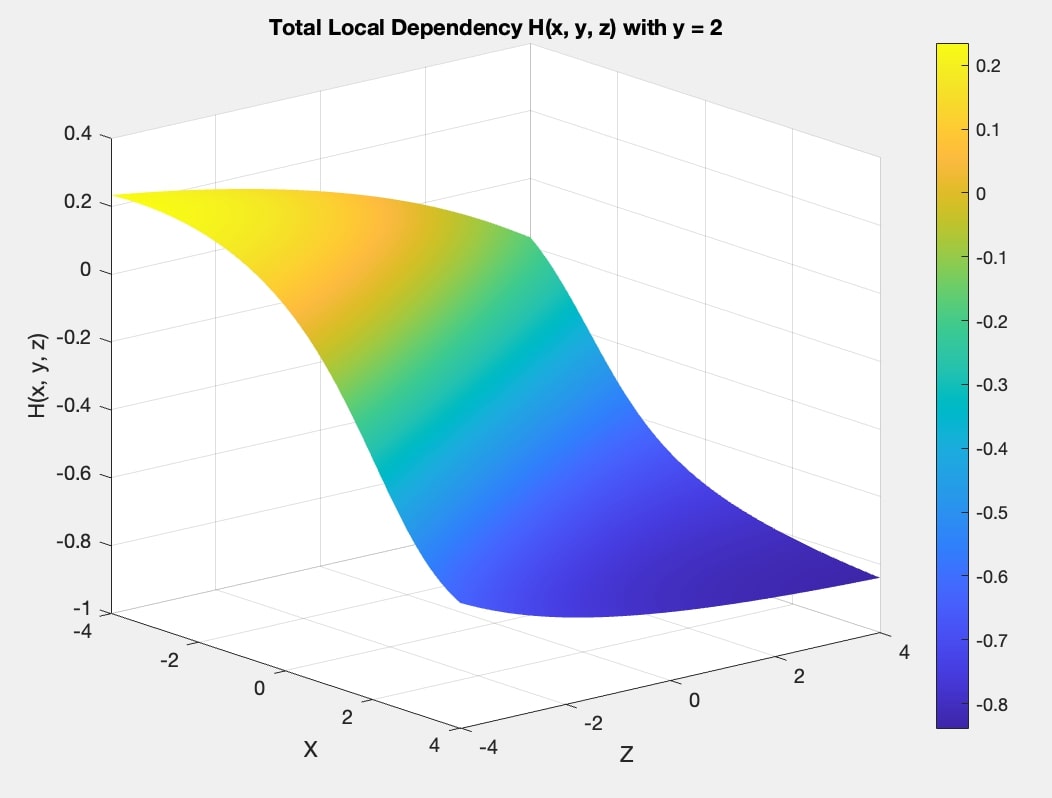}
    \end{subfigure}
    \caption{Graph of \( H(x,y,z) \) for \( y = 2 \),  \( x,z\in (-4,4)\).}
    \label{fig:8}
\end{figure}

\begin{figure}[H]
    \centering
    \begin{subfigure}[t]{0.45\textwidth}
        \centering
        \includegraphics[width=\textwidth]{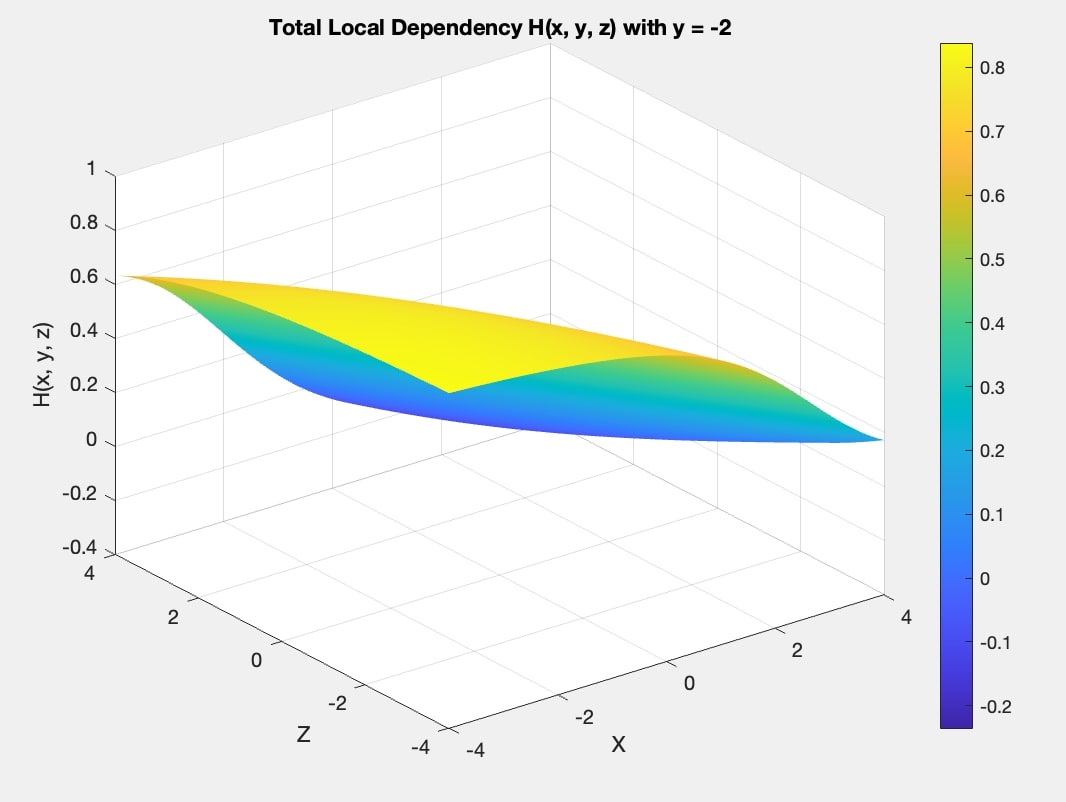}
    \end{subfigure}
    \hfill
    \begin{subfigure}[t]{0.45\textwidth}
        \centering
        \includegraphics[width=\textwidth]{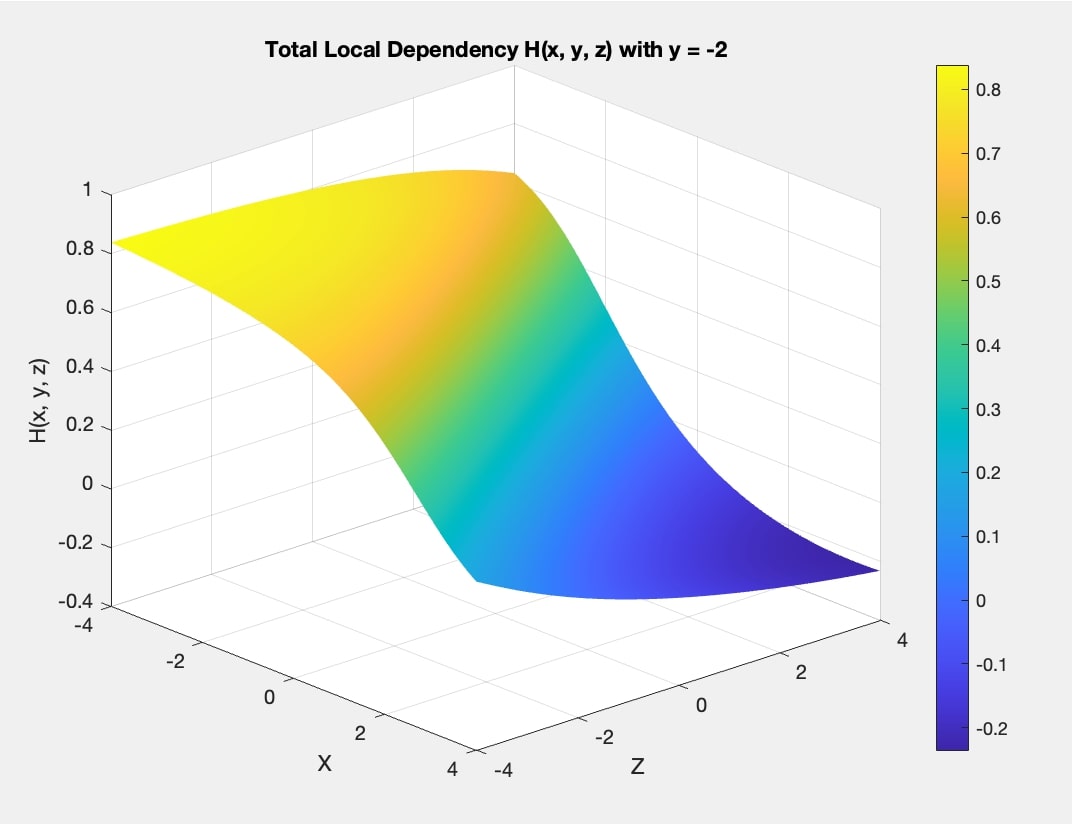}
    \end{subfigure}
    \caption{Graph of \( H(x,y,z) \) for \( y = -2 \),  \( x,z\in (-4,4)\).}
    \label{fig:9}
\end{figure}

\begin{figure}[H]
    \centering
    \begin{subfigure}[t]{0.45\textwidth}
        \centering
        \includegraphics[width=\textwidth]{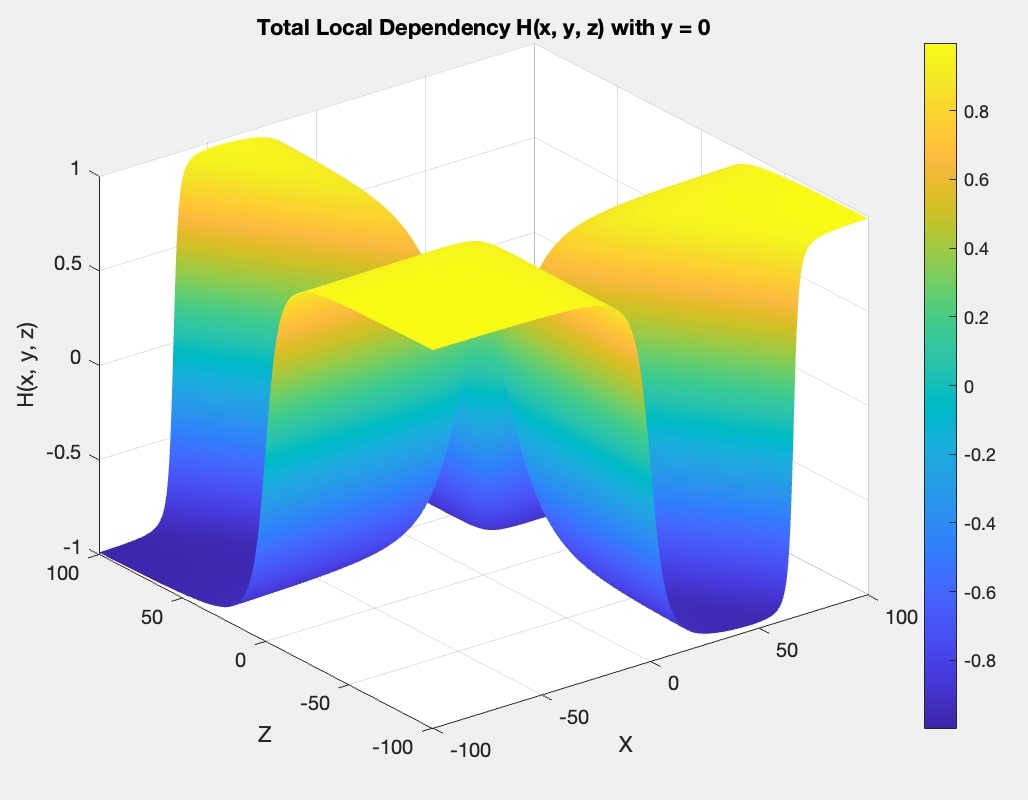}
    \end{subfigure}
    \hfill
    \begin{subfigure}[t]{0.45\textwidth}
        \centering
        \includegraphics[width=\textwidth]{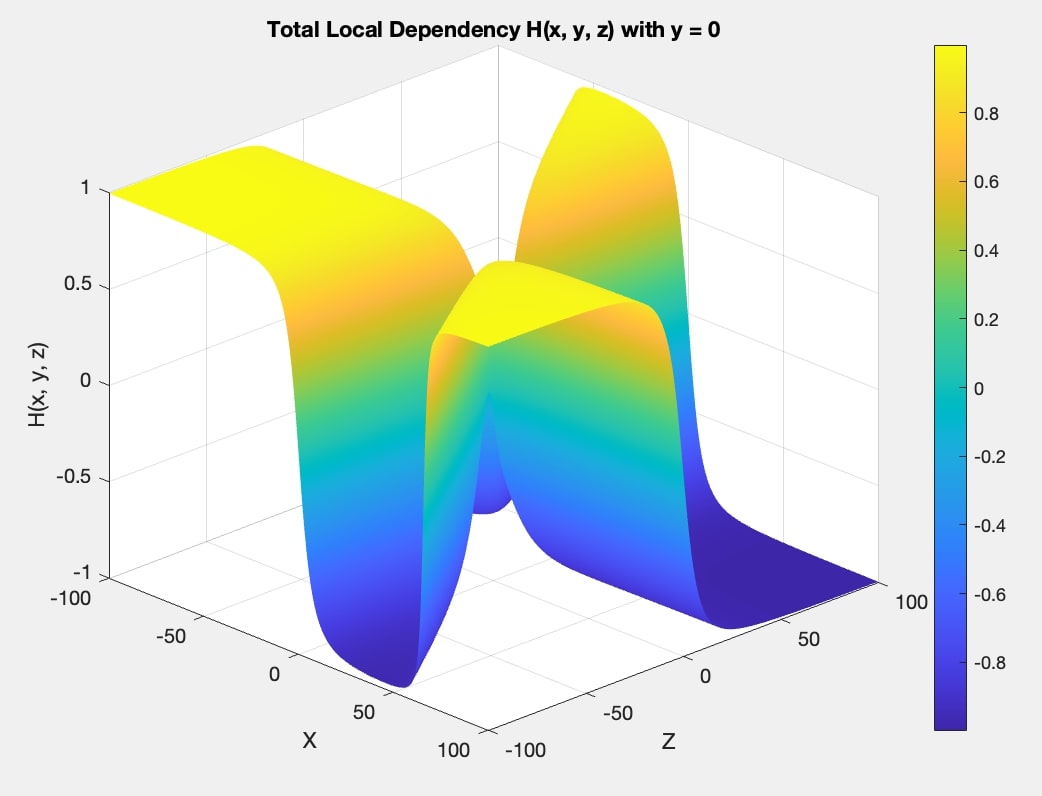}
    \end{subfigure}
    \caption{Graph of \( H(x,y,z) \) for \( y = 0 \),  \( x,z\in (-100,100)\).}
    \label{fig:10}
\end{figure}

Similar graphs are observed for total local dependence $H(x,y,z)$ for different fixed values of z and varying intervals for $x$ and $y$.

To provide a different insight to the total local dependence measure, the covariance between pairs of variables were slightly increased and the covariance matrix for these computations was taken as follows:

Covariance Matrix ($\Sigma$): 
\begin{equation*}
\Sigma = 
\begin{bmatrix}
\sigma_{X}^{2} & Cov(X,Y) & Cov(X,Z) \\ 
Cov(Y,X) & \sigma_{y}^{2} & Cov(Y,Z) \\ 
Cov(Z,X) & Cov(Z,Y) & \sigma_{z}^{2}
\end{bmatrix}
= 
\begin{bmatrix}
1.0 & 0.8 & 0.6 \\ 
0.8 & 1.0 & 0.4 \\ 
0.6 & 0.4 & 1.0
\end{bmatrix}%
\end{equation*}
This matrix defines the variance for each variable ($\sigma_x^2 = \sigma_y^2= \sigma_z^2 = 1$) and the covariance between pairs of variables. The covariance values are taken again arbitrarily.

To visualize the total local dependence in a 3D setting for the new covariance matrix, we again fixed $x$, $y$ and $z$ at $0$ in separate cases and varying the others over a specific interval.

\begin{figure}[H]
    \centering
    \begin{subfigure}[t]{0.45\textwidth}
        \centering
        \includegraphics[width=\textwidth]{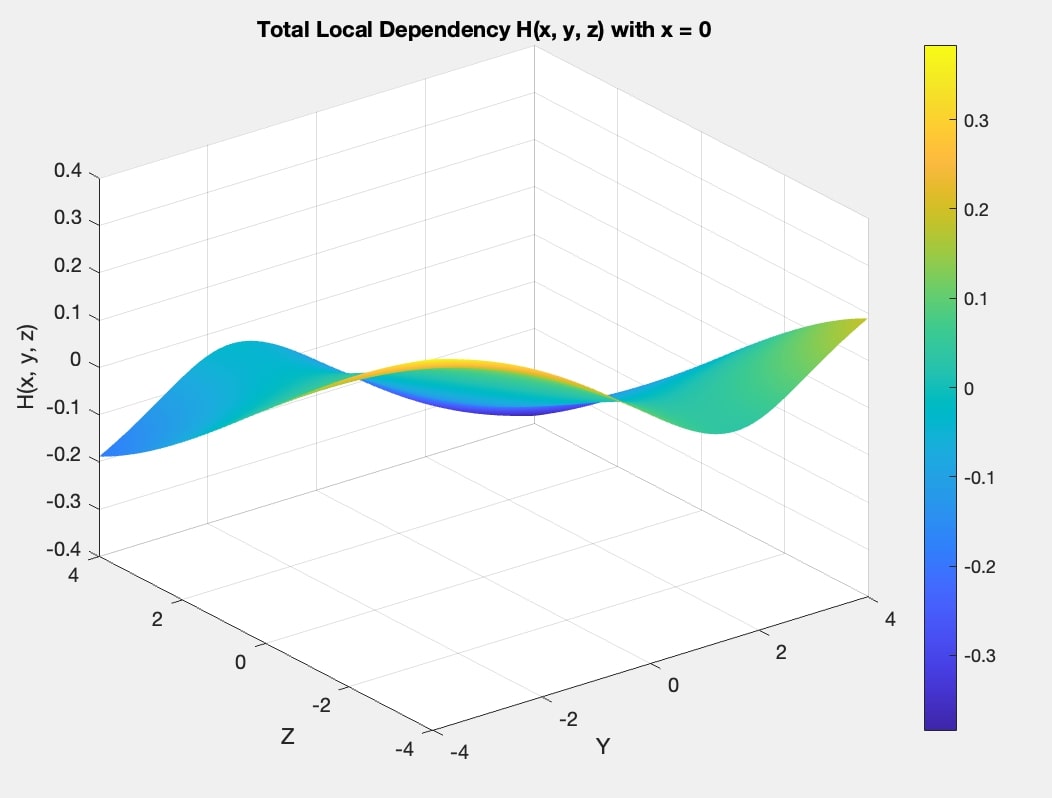}
    \end{subfigure}
    \hfill
    \begin{subfigure}[t]{0.45\textwidth}
        \centering
        \includegraphics[width=\textwidth]{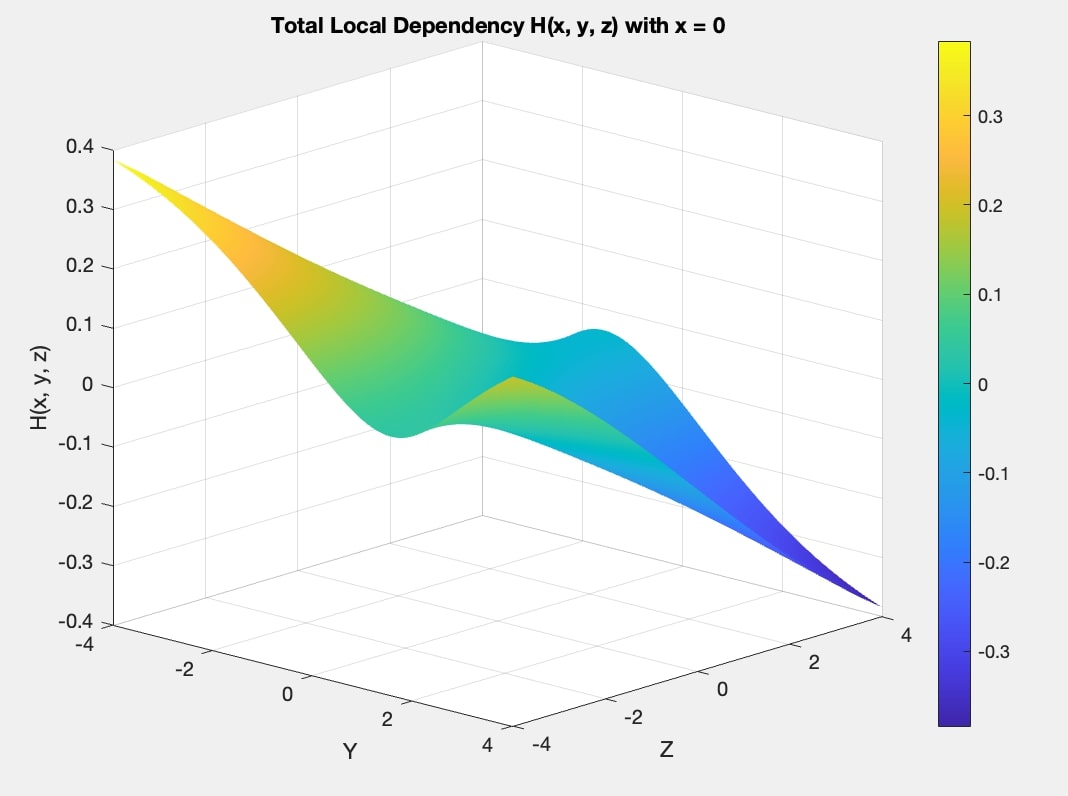}
    \end{subfigure}
    \caption{Total local dependence dependence with increased sigma values for \( x = 0 \),  \( y,z\in (-4,4)\).}
    \label{fig:11}
\end{figure}

\begin{figure}[H]
    \centering
    \begin{subfigure}[t]{0.45\textwidth}
        \centering
        \includegraphics[width=\textwidth]{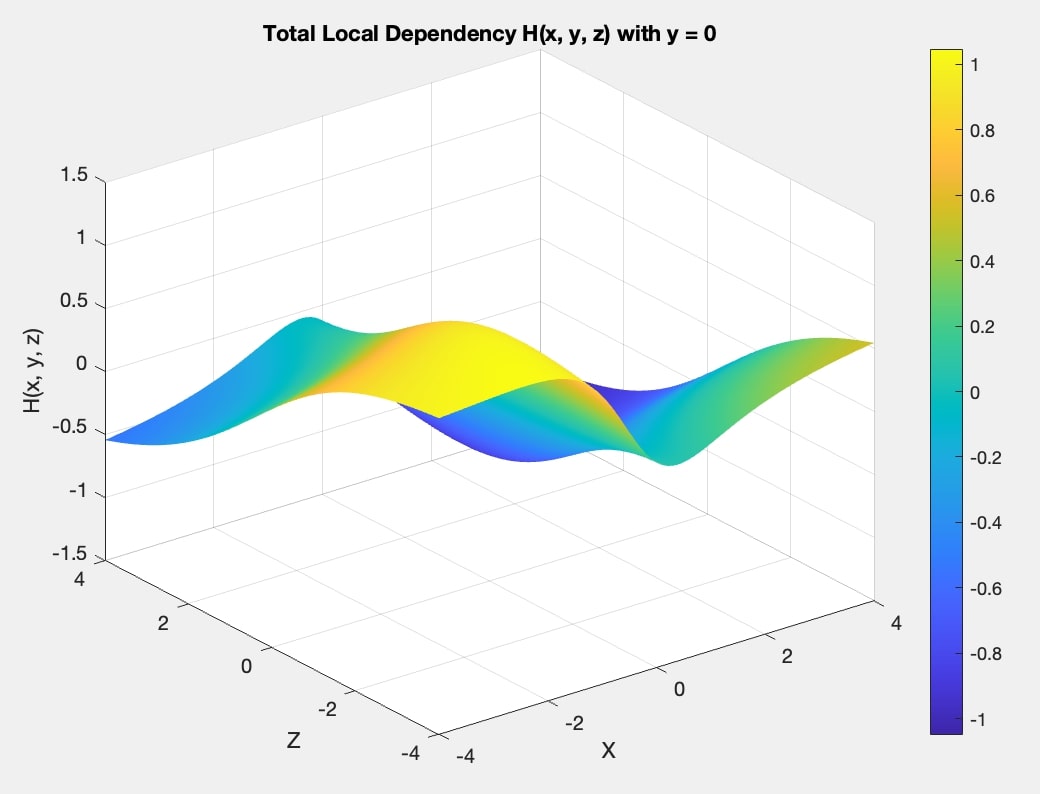}
    \end{subfigure}
    \hfill
    \begin{subfigure}[t]{0.45\textwidth}
        \centering
        \includegraphics[width=\textwidth]{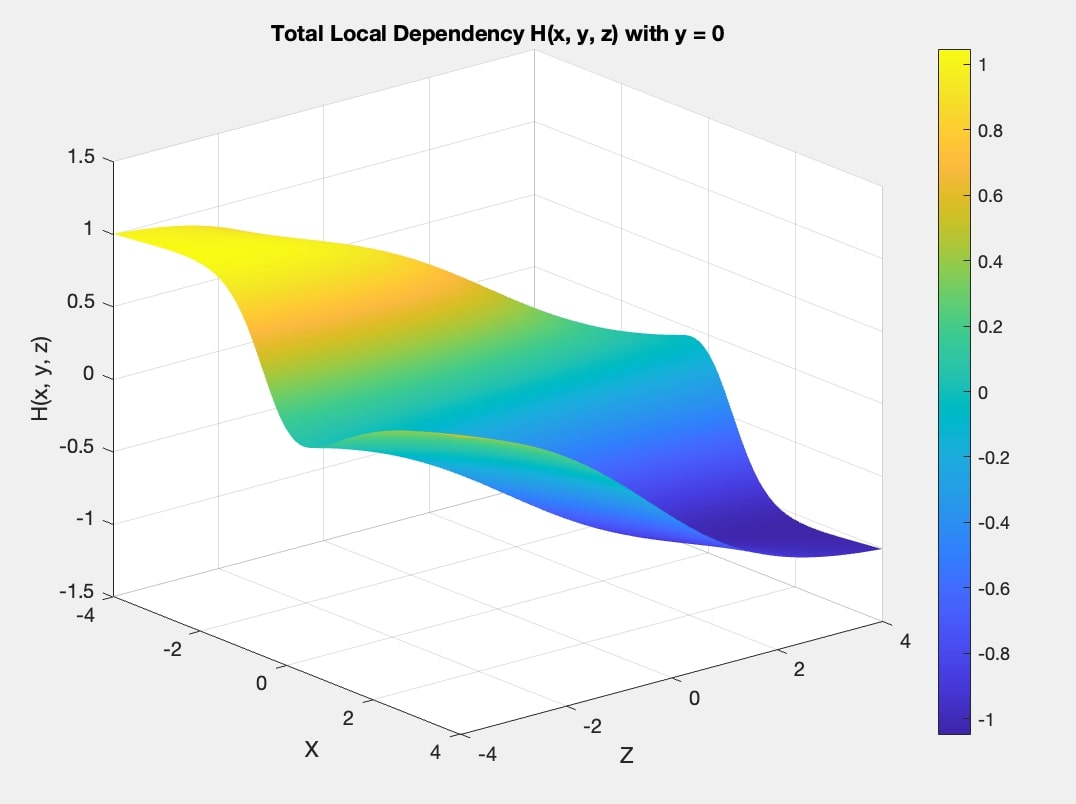}
    \end{subfigure}
    \caption{Total local dependence dependence with increased sigma values for \( y = 0 \),  \( x,z\in (-4,4)\).}
    \label{fig:12}
\end{figure}

\begin{figure}[H]
    \centering
    \begin{subfigure}[t]{0.45\textwidth}
        \centering
        \includegraphics[width=\textwidth]{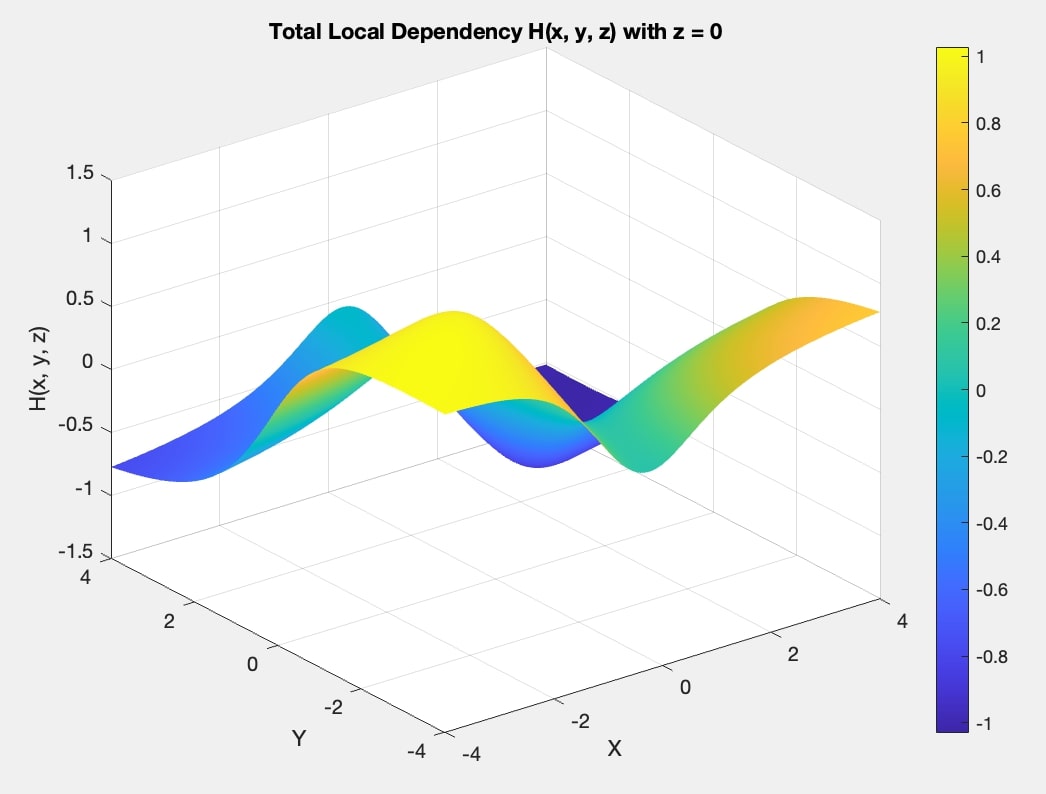}
    \end{subfigure}
    \hfill
    \begin{subfigure}[t]{0.45\textwidth}
        \centering
        \includegraphics[width=\textwidth]{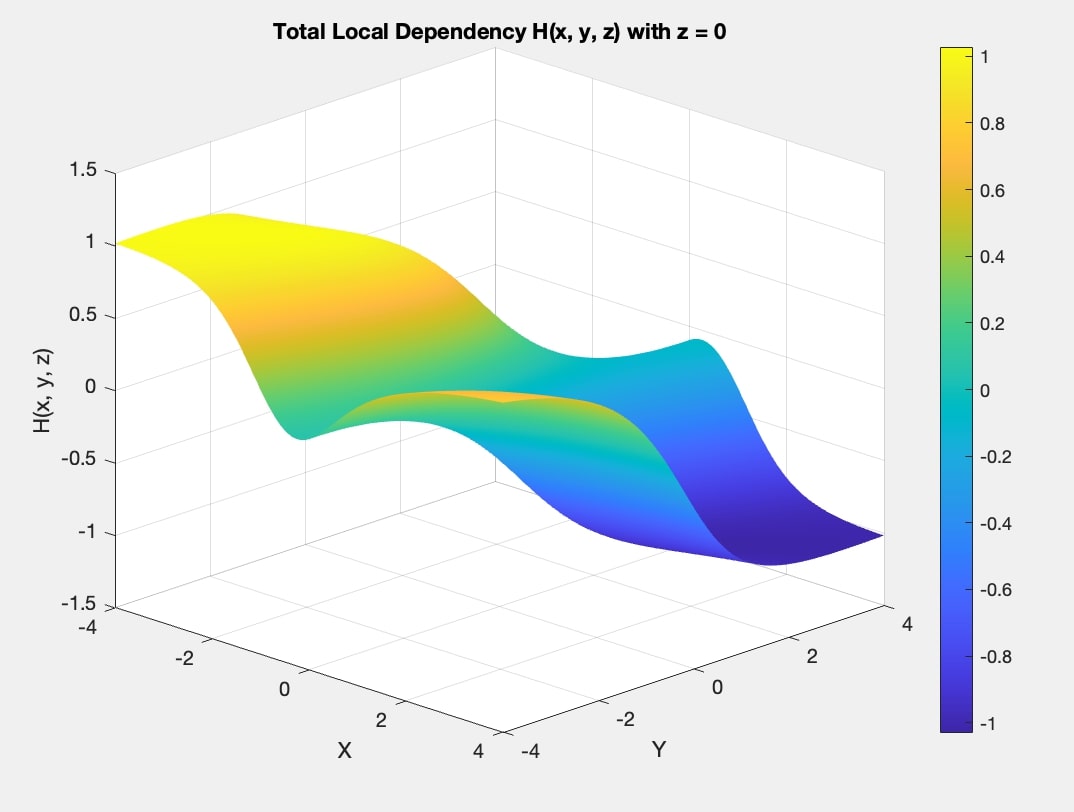}
    \end{subfigure}
    \caption{Total local dependence dependence with increased sigma values for \( z = 0 \),  \( x,y\in (-4,4)\).}
    \label{fig:13}
\end{figure}

\subsection{Tables}

We provide the calculated values of the local dependence function in some particular points for three-variate normal distribution with the vector of means (\ref{nn1}) and the covariance matrix (\ref{nn2}).

\begin{table}[H]
    \centering
    \caption{Numerical values of H(x, y, z) for the given $\mu$ and $\Sigma$} \label{tab:table1}
    \begin{tabular}{rrrr|rrrr}
        \toprule
        \textbf{x} & \textbf{y} & \textbf{z} & \textbf{H(x, y, z)} & \textbf{x} & \textbf{y} & \textbf{z} & \textbf{H(x, y, z)} \\
        \midrule
        0   & 0   & 0    & 0.0000   &  0   & 100 & 0   & -0.0344 \\
        0   & 0   & 1    & -0.1245  &  0   & 100 & 1   & -0.1520 \\
        0   & 0   & -1   & 0.1245   &  0   & 100 & -1  & 0.0837  \\
        0   & 0   & 10   & -0.2861  &  0   & 100 & 10  & -0.7855 \\
        0   & 0   & -10  & 0.2861   &  0   & 100 & -10 & 0.7413  \\
        0   & 0   & 20   & -0.1803  &  0   & 100 & 20  & -0.9325 \\
        0   & 0   & -20  & 0.1803   &  0   & 100 & -20 & 0.9058  \\
        0   & 0   & 100  & -0.0396  &  0   & 100 & 100 & -0.9965 \\
        0   & 0   & -100 & 0.0396   &  0   & 100 & -100 & 0.9886  \\ \midrule
        0   & 1   & 0    & -0.3005  &  1   & 0   & 0   & -0.3638 \\
        0   & 1   & 1    & -0.4209  &  1   & 0   & 1   & -0.1756 \\
        0   & 1   & -1   & -0.1623  &  1   & 0   & -1  & -0.2695 \\
        0   & 1   & 10   & -0.5319  &  1   & 0   & 10  & -0.3831 \\
        0   & 1   & -10  & 0.0723   &  1   & 0   & -10 & 0.1857  \\
        0   & 1   & 20   & -0.4601  &  1   & 0   & 20  & -0.2989 \\
        0   & 1   & -20  & -0.1000  &  1   & 0   & -20 & 0.0591  \\
        0   & 1   & 100  & -0.3522  &  1   & 0   & 100 & -0.1708 \\
        0   & 1   & -100 & -0.2768  &  1   & 0   & -100 & -0.0923 \\ \midrule
        0   & -1  & 0    & 0.3005   &  -1  & 0   & 0   & 0.1756  \\
        0   & -1  & 1    & 0.1623   &  -1  & 0   & 1   & 0.0581  \\
        0   & -1  & -1   & 0.4209   &  -1  & 0   & -1  & 0.2695  \\
        0   & -1  & 10   & -0.0723  &  -1  & 0   & 10  & -0.1857 \\
        0   & -1  & -10  & 0.5319   &  -1  & 0   & -10 & 0.3831  \\
        0   & -1  & 20   & -0.1000  &  -1  & 0   & 20  & -0.0591 \\
        0   & -1  & -20  & -0.4601  &  -1  & 0   & -20 & 0.2989  \\
        0   & -1  & 100  & -0.2768  &  -1  & 0   & 100 & -0.0923 \\
        0   & -1  & -100 & -0.3522  &  -1  & 0   & -100 & 0.1708  \\
        \bottomrule
    \end{tabular}
\end{table}

\begin{table}[H]
    \centering
    \caption{Numerical values of H(x, y, z) for the given $\mu$ and $\Sigma$}\label{tab:table2}
    \begin{tabular}{rrrr|rrrr}
        \toprule
        \textbf{x} & \textbf{y} & \textbf{z} & \textbf{H(x, y, z)} & \textbf{x} & \textbf{y} & \textbf{z} & \textbf{H(x, y, z)} \\
        \midrule
        10  & 0   & 0    & -0.2682  &  -10 & 0   & 0   & 0.2682  \\
        10  & 0   & 1    & -0.3600  &  -10 & 0   & 1   & 0.1731  \\
        10  & 0   & -1   & -0.1731  &  -10 & 0   & -1  & 0.3600  \\
        10  & 0   & 10   & -0.7839  &  -10 & 0   & 10  & -0.3648 \\
        10  & 0   & -10  & 0.3648   &  -10 & 0   & -10 & 0.7839  \\
        10  & 0   & 20   & -0.8417  &  -10 & 0   & 20  & -0.4921 \\
        10  & 0   & -20  & -0.4921  &  -10 & 0   & -20 & 0.8417  \\
        10  & 0   & 100  & -0.8205  &  -10 & 0   & 100 & -0.7726 \\
        10  & 0   & -100 & -0.7726  &  -10 & 0   & -100 & 0.8205  \\ \midrule
        10  & 1   & 0    & -0.5912  &  -10 & 1   & 0   & -0.0476 \\
        10  & 1   & 1    & -0.6498  &  -10 & 1   & 1   & -0.1200 \\
        10  & 1   & -1   & -0.5217  &  -10 & 1   & -1  & 0.0321  \\
        10  & 1   & 10   & -0.8614  &  -10 & 1   & 10  & -0.4376 \\
        10  & 1   & -10  & 0.1877   &  -10 & 1   & -10 & 0.6492  \\
        10  & 1   & 20   & -0.8901  &  -10 & 1   & 20  & -0.3850 \\
        10  & 1   & -20  & -0.5588  &  -10 & 1   & -20 & 0.7641  \\
        10  & 1   & 100  & -0.8743  &  -10 & 1   & 100 & -0.6762 \\
        10  & 1   & -100 & -0.8328  &  -10 & 1   & -100 & 0.7324  \\ \midrule
        10  & -1  & 0    & 0.0476   &  -10 & -1  & 0   & 0.5912  \\
        10  & -1  & 1    & -0.0321  &  -10 & -1  & 1   & 0.5217  \\
        10  & -1  & -1   & 0.1200   &  -10 & -1  & -1  & 0.6498  \\
        10  & -1  & 10   & -0.6492  &  -10 & -1  & 10  & -0.1877 \\
        10  & -1  & -10  & 0.4376   &  -10 & -1  & -10 & 0.8614  \\
        10  & -1  & 20   & -0.7641  &  -10 & -1  & 20  & -0.5588 \\
        10  & -1  & -20  & -0.3850  &  -10 & -1  & -20 & 0.8901  \\
        10  & -1  & 100  & -0.7324  &  -10 & -1  & 100 & -0.8328 \\
        10  & -1  & -100 & -0.6762  &  -10 & -1  & -100 & 0.8743  \\ \midrule
        10  & 10  & 0    & -0.9853  &  -10 & 10  & 0   & 0.8143  \\
        10  & 10  & 1    & -0.9861  &  -10 & 10  & 1   & 0.8076  \\
        10  & 10  & -1   & -0.9843  &  -10 & 10  & -1  & 0.8200  \\
        10  & 10  & 10   & -0.9878  &  -10 & 10  & 10  & 0.5955  \\
        10  & 10  & -10  & -0.9087  &  -10 & 10  & -10 & 0.8432  \\
        10  & 10  & 20   & -0.9873  &  -10 & 10  & 20  & -0.8265 \\
        10  & 10  & -20  & 0.5564   &  -10 & 10  & -20 & 0.8664  \\
        10  & 10  & 100  & -0.9829  &  -10 & 10  & 100 & -0.9089 \\
        10  & 10  & -100 & -0.9589  &  -10 & 10  & -100 & 0.8642  \\
        \bottomrule
    \end{tabular}
\end{table}

\bigskip

For comparison of values of the joint p.d.f $f(x,y,z)$ and $H(x,y,z)$ we provide calculated values at some particular points.

\begin{table}[H]
    \centering 
    \caption{Numerical values of $f(x, y, z)$ and $H(x, y, z)$ for the given $\mu$ and $\Sigma$} \label{tab:table3}
    \begin{tabular}{rrrrr|rrrrr}
        \toprule
        \textbf{x} & \textbf{y} & \textbf{z} & \textbf{f(x, y, z)} & \textbf{H(x, y, z)} & \textbf{x} & \textbf{y} & \textbf{z} & \textbf{f(x, y, z)} & \textbf{H(x, y, z)} \\
        \midrule
        -1.00 & 0.80  & 1.00  & 0.0101 & -0.1831 & -0.80 & -1.00 & -1.00 & 0.0381 & 0.5324 \\
        -0.80 & -1.00 & -0.80 & 0.0426 & 0.5188 & -0.80 & -0.80 & -1.00 & 0.0415 & 0.4839 \\
        -0.80 & -1.00 & -0.50 & 0.0460 & 0.4955 & -0.80 & -0.80 & -0.80 & 0.0471 & 0.4697 \\
        -0.80 & -1.00 & 0.00  & 0.0410 & 0.4490 & -0.80 & -0.80 & -0.50 & 0.0521 & 0.4458 \\
        -0.80 & -1.00 & 0.50  & 0.0270 & 0.3930 & -0.80 & -0.80 & 0.00  & 0.0484 & 0.3984 \\
        -0.80 & -1.00 & 0.80  & 0.0182 & 0.3555 & -0.80 & -0.80 & 0.50  & 0.0332 & 0.3421 \\
        -0.80 & -1.00 & 1.00  & 0.0131 & 0.3293 & -0.80 & -0.80 & 0.80  & 0.0229 & 0.3049 \\
        -0.80 & -0.80 & -1.00 & 0.0415 & 0.4839 & -0.80 & -0.50 & -1.00 & 0.0423 & 0.4010 \\
        -0.80 & -0.80 & -0.80 & 0.0471 & 0.4697 & -0.80 & -0.50 & -0.80 & 0.0492 & 0.3857 \\
        -0.80 & -0.80 & -0.50 & 0.0521 & 0.4458 & -0.80 & -0.50 & -0.50 & 0.0564 & 0.3601 \\
        -0.80 & -0.80 & 0.00  & 0.0484 & 0.3984 & -0.80 & -0.50 & 0.00  & 0.0556 & 0.3106 \\
        -0.80 & -0.80 & 0.50  & 0.0332 & 0.3421 & -0.80 & -0.50 & 0.50  & 0.0405 & 0.2533 \\
        -0.80 & -0.80 & 0.80  & 0.0229 & 0.3049 & -0.80 & -0.50 & 0.80  & 0.0290 & 0.2161 \\
        -0.80 & -0.80 & 1.00  & 0.0168 & 0.2790 & -0.80 & -0.50 & 1.00  & 0.0218 & 0.1906 \\
        -0.80 & -0.50 & -1.00 & 0.0423 & 0.4010 & -0.80 & 0.00  & -1.00 & 0.0325 & 0.2457 \\
        -0.80 & -0.50 & -0.80 & 0.0492 & 0.3857 & -0.80 & 0.00  & -0.80 & 0.0394 & 0.2279 \\
        -0.80 & -0.50 & -0.50 & 0.0564 & 0.3601 & -0.80 & 0.00  & -0.50 & 0.0479 & 0.1989 \\
        -0.80 & -0.50 & 0.00  & 0.0556 & 0.3106 & -0.80 & 0.00  & 0.00  & 0.0523 & 0.1448 \\
        -0.80 & -0.50 & 0.50  & 0.0405 & 0.2533 & -0.80 & 0.00  & 0.50  & 0.0421 & 0.0850 \\
        -0.80 & -0.50 & 0.80  & 0.0290 & 0.2161 & -0.80 & 0.00  & 0.80  & 0.0320 & 0.0476 \\
        \bottomrule
    \end{tabular}
\end{table}

\begin{table}[H]
    \centering
    \caption{Numerical values of $f(x, y, z)$ and $H(x, y, z)$ for the given $\mu$ and $\Sigma$} \label{tab:table4}
    \begin{tabular}{rrrrr|rrrrr}
        \toprule
        \textbf{x} & \textbf{y} & \textbf{z} & \textbf{f(x, y, z)} & \textbf{H(x, y, z)} & \textbf{x} & \textbf{y} & \textbf{z} & \textbf{f(x, y, z)} & \textbf{H(x, y, z)} \\
        \midrule
        0.00  & -0.80 & 1.00  & 0.0199 & 0.1179 & 0.00  & -0.50 & -1.00 & 0.0448 & 0.2874 \\
        0.00  & -0.50 & -0.80 & 0.0536 & 0.2656 & 0.00  & -0.50 & -0.50 & 0.0638 & 0.2306 \\
        0.00  & -0.50 & 0.00  & 0.0671 & 0.1672 & 0.00  & -0.50 & 0.50  & 0.0522 & 0.1012 \\
        0.00  & -0.50 & 0.80  & 0.0388 & 0.0621 & 0.00  & -0.50 & 1.00  & 0.0300 & 0.0368 \\
        0.00  & 0.00  & -1.00 & 0.0440 & 0.1245 & 0.00  & 0.00  & -0.80 & 0.0548 & 0.1011 \\
        0.00  & 0.00  & -0.50 & 0.0693 & 0.0643 & 0.00  & 0.00  & 0.00  & 0.0806 & 0.0000 \\
        0.00  & 0.00  & 0.50  & 0.0693 & -0.0643 & 0.00  & 0.00  & 0.80  & 0.0548 & -0.1011 \\
        0.00  & 0.00  & 1.00  & 0.0440 & -0.1245 & 0.00  & 0.50  & -1.00 & 0.0300 & -0.0368 \\
        0.00  & 0.50  & -0.80 & 0.0388 & -0.0621 & 0.00  & 0.50  & -0.50 & 0.0522 & -0.1012 \\
        0.00  & 0.50  & 0.00  & 0.0671 & -0.1672 & 0.00  & 0.50  & 0.50  & 0.0638 & -0.2306 \\
        0.00  & 0.50  & 0.80  & 0.0536 & -0.2656 & 0.00  & 0.50  & 1.00  & 0.0448 & -0.2874 \\
        0.00  & 0.80  & -1.00 & 0.0199 & -0.1179 & 0.00  & 0.80  & -0.80 & 0.0264 & -0.1445 \\
        0.00  & 0.80  & -0.50 & 0.0369 & -0.1851 & 0.00  & 0.80  & 0.00  & 0.0504 & -0.2527 \\
        0.00  & 0.80  & 0.50  & 0.0509 & -0.3163 & 0.00  & 0.80  & 0.80  & 0.0443 & -0.3510 \\
        0.00  & 0.80  & 1.00  & 0.0380 & -0.3723 & 0.00  & 1.00  & -1.00 & 0.0316 & -0.4209 \\
        0.50  & -1.00 & -1.00 & 0.0181 & 0.3279 & 0.50  & -0.80 & -1.00 & 0.0232 & 0.2808 \\
        0.50  & -1.00 & -0.80 & 0.0211 & 0.3028 & 0.50  & -0.80 & -0.80 & 0.0274 & 0.2561 \\
        0.50  & -1.00 & -0.50 & 0.0243 & 0.2633 & 0.50  & -0.80 & -0.50 & 0.0323 & 0.2173 \\
        0.50  & -1.00 & 0.00  & 0.0241 & 0.1952 & 0.50  & -0.80 & 0.00  & 0.0333 & 0.1507 \\
        0.50  & -1.00 & 0.50  & 0.0176 & 0.1278 & 0.50  & -0.80 & 0.50  & 0.0254 & 0.0853 \\
        \bottomrule
    \end{tabular}
\end{table}

\bigskip

\section{The multivariate extension of local dependence function}

\bigskip

$X_{1},X_{2},....,X_{n}$ be dependent random variables with $EX_{1}=\mu
_{1},EX_{2}=\mu _{2},...,EX_{n}=\mu _{n}$ and $Var(X_{1})=\sigma
_{1}^{2},...,Var(X_{n})=\sigma _{n}^{2}.$ The local dependence function is

\begin{eqnarray*}
H(x_{1},x_{2},...,x_{n}) &=&[E((X_{1}-E(X_{1}\mid X_{2}=x_{2},...,X_{n}=x_{n})) \\
\times E((X_{2}-E(X_{2} &\mid &X_{1}=x_{1},X_{3}=x_{3},...,X_{n}=x_{n})) \\
\cdots \times E((X_{n}-E(X_{n} &\mid &X_{1}=x_{1},X_{2}=x_{2},...,X_{n-1}=x_{n-1}))]/ \\
&&(\sqrt{E[(X_{1}-E(X_{1}\mid X_{2}=x_{2},...,X_{n}=x_{n}))^{2}]} \\
&&\times \sqrt{E[(X_{2}-E(X_{2}\mid X_{1}=x_{1},X_{3}=x_{3},...,X_{n}=x_{n}))^{2}]}\cdots \\
&&\times \sqrt{E[(X_{n}-E(X_{n}\mid X_{1}=x_{1},X_{2}=x_{2},...,X_{n-1}=x_{n-1}))^{2}]})
\end{eqnarray*}

The function $H(x_{1},x_{2},...,x_{n})$ can be expressed as 
\begin{eqnarray*}
&&H(x_{1},x_{2},...,x_{n}) \\
&=&[E((X_{1}-EX_{1}+(EX_{1}-E(X_{1}\mid X_{2}=x_{2},...,X_{n}=x_{n}))) \\
&&\times E((X_{2}-EX_{2} \\
+(EX_{2}-E(X_{2} &\mid &X_{1}=x_{1},X_{3}=x_{3},...,X_{n}=x_{n}))) \\
&&\cdots \times E((X_{n}-EX_{n} \\
+(EX_{n}-E(X_{n} &\mid &X_{1}=x_{1},X_{2}=x_{2},...,X_{n-1}=x_{n-1})))]/
\end{eqnarray*}
\begin{equation*}
(\sqrt{E[(X_{1}-EX_{1}+(EX_{1}-E(X_{1}\mid
X_{2}=x_{2},...,X_{n}=x_{n})))^{2}]}
\end{equation*}
\begin{eqnarray*}
&&\times \sqrt{E[(X_{2}-EX_{2}+(EX_{2}-E(X_{2}\mid X_{1}=x_{1},X_{3}=x_{3},...,X_{n}=x_{n})))^{2}]}\cdots \\
&&\times \sqrt{E[(X_{n}-EX_{n}+(EX_{n}-E(X_{n}\mid X_{1}=x_{1},X_{2}=x_{2},...,X_{n-1}=x_{n-1})))^{2}]})
\end{eqnarray*}

In the case of $n=4,$ we can write
\begin{eqnarray*}
&&H(x_{1},x_{2},x_{3},x_{4}) \\
&=&[\rho _{X_{1},X_{2},X_{3},X_{4}}+\rho _{X_{1},X_{2},X_{3}}\varphi_{X_{4}}(x_{1},x_{2},x_{3}) \\
&&+\rho _{X_{2},X_{3},X_{4}}\varphi _{X_{1}}(x_{2},x_{3},x_{4})+\rho_{X_{1},X_{3},X_{4}}\varphi _{X_{2}}(x_{1},x_{3},x_{4}) \\
&&\rho _{X_{1},X_{2},X_{4}}\varphi _{X_{3}}(x_{1},x_{2},x_{4})+\rho_{X_{1},X_{2}}\varphi _{X_{3}}(x_{1},x_{2},x_{4})\varphi_{X_{4}}(x_{1},x_{2},x_{3})
\end{eqnarray*}
\begin{eqnarray*}
&&+\rho _{X_{1},X_{3}}\varphi _{X_{2}}(x_{1},x_{2},x_{4})\varphi_{X_{4}}(x_{1},x_{2},x_{3}) \\
&&+\rho _{X_{1},X_{4}}\varphi _{X_{3}}(x_{1},x_{2},x_{4})\varphi_{X_{2}}(x_{1},x_{3},x_{4}) \\
&&+\rho _{X_{2},X_{3}}\varphi _{X_{1}}(x_{2},x_{3},x_{4})\varphi_{X_{4}}(x_{1},x_{2},x_{3}) \\
&&+\rho _{X_{2},X_{4}}\varphi _{X_{3}}(x_{1},x_{2},x_{4})\varphi_{X_{1}}(x_{2},x_{3},x_{4}) \\
&&+\rho _{X_{3},X_{4}}\varphi _{X_{1}}(x_{2},x_{3},x_{4})\varphi_{X_{2}}(x_{1},x_{3},x_{4}) \\
&&+\varphi _{X_{1}}(x_{2},x_{3},x_{4})\varphi_{X_{2}}(x_{1},x_{3},x_{4})\varphi _{X_{3}}(x_{1},x_{2},x_{4})\varphi_{X_{4}}(x_{1},x_{2},x_{3})]/ \\
&&[\sqrt{1+\varphi _{X_{1}}^{2}(x_{2},x_{3},x_{4})}\sqrt{1+\varphi_{X_{2}}^{2}(x_{1},x_{3},x_{4})}\times \\
&&\times \sqrt{1+\varphi _{X_{3}}^{2}(x_{1},x_{2},x_{4})}\sqrt{1+\varphi_{X_{4}}^{2}(x_{1},x_{2},x_{3})}],
\end{eqnarray*}

where
\begin{equation*}
\rho _{X_{1},X_{2},X_{3},X_{4}}=\frac{E(X_{1}-\mu _{1})(X_{2}-\mu_{2})(X_{3}-\mu _{3})(X_{4}-\mu _{4})}{\sigma _{1}\sigma_{2}\sigma_{3}\sigma _{4}}
\end{equation*}
\begin{eqnarray*}
\rho _{X_{i},X_{j},X_{k}} &=&\frac{E(X_{i}-\mu _{i})(X_{j}-\mu
_{j})(X_{k}-\mu _{k})}{\sigma _{i}\sigma _{j}\sigma _{k}} \\
\rho _{X_{i},X_{j}} &=&\frac{E(X_{i}-\mu _{i})(X_{j}-\mu _{j})}{\sigma
_{i}\sigma _{j}}
\end{eqnarray*}

\subsection{Properties of multivariate local dependence function}

The local dependence function $H(x_{1},x_{2},...,x_{n})$ possess similar properties to that $H(x,y,z)$ has.

\begin{theorem}
1. \ If $X_{1},X_{2},...,X_{n}$ are jointly independent
random variables, i.e. the connecting copula for a joit distribution function is $C(t_{1},t_{2},...,t_{n})=\prod{i=1}^{n}t_{i},$ $
0\leq t_{i}\leq 1,i=1,2...,n.$ Then $H(x_{1},x_{2},...,x_{n})=0$ for all $(x_{1},x_{2},...,x_{n})\in \mathbb{R}^{n}.$ \ \newline
2. $\left\vert H(x_{1},x_{2},...,x_{n})\right\vert \leq 1,$ for all $(x_{1},x_{2},...,x_{n})\in \mathbb{R}^{n}$ \newline
3. \ Let the random variables $X_{i_{1}},X_{i_{2}},...,X_{i_{n-1}}$ ,$1\leq
i_{1}<i_{2}<\cdots <i_{k}\leq n,2\leq k\leq n-1$ be independent, i.e. any two, three etc. $n-1$ of random variables $X_{1},X_{2},...,X_{n}$ are independent. Then
\begin{equation*}
	H(x_{1},x_{2},...,x_{n})=\frac{\rho
		_{X_{1},X_{2},...,X_{n}}+\prod_{i=1}^{n}\varphi
		_{X_{i}}(x_{1},x_{2},...,x_{i-1},x_{i+1},...,x_{n})}{\prod_{i=1}^{n}%
		\sqrt{1+\varphi _{X_{i}}^{2}(x_{1},x_{2},...,x_{i-1},x_{i+1},...,x_{n})}}
\end{equation*}
for all $(x,y,z)\in \mathbb{R}^{3},$ where \newline
\begin{eqnarray*}
\rho _{X_{1},X_{2},...,X_{n}} &=&E(X_{1}-EX_{1})\cdots (X_{n}-EX_{n}) \\
&&\varphi _{X_{i}}(x_{1},x_{2},...,x_{i-1},x_{i+1},...,x_{n}) \\	&=&\frac{EX_{i}-E(X_{i}|X_{1}=x_{1},X_{2}=x_{2},...,X_{i-1}=x_{i-1},X_{i+1}=x_{i+1},...,X_{n}=x_{n})}{\sigma _{i}} \\
	i &=&1,2,...,n.
\end{eqnarray*}
4. \ Let $(x_{1},x_{2},...,x_{n})\in \mathbb{R}^{3}$ satisfy the system of equations
\end{theorem}

\begin{equation}
\left\{
\begin{array}{c}
	EX_{1}=E(X_{1}\mid X_{2}=x_{2},...,X_{n}=x_{n}) \\
	EX_{2}=E(X_{2}\mid X_{1}=x_{1},X_{3}=x_{3},...,X_{n}=x_{n}) \\
	\cdots \\
	EX_{n}=E(X_{n}\mid X_{1}=x_{1},X_{2}=x_{2},...,X_{n-1}=x_{n-1})%
\end{array}
\right.
\end{equation}
\ then $H(x_{1},x_{2},...,x_{n})=\rho _{X_{1},X_{2},...,X_{n}}.$\newline

\bigskip

\begin{conclusion}
The local dependence function is important in many applications of statistics and probability theory. The classical dependence measures are scalars and can not describe the dependence between random variables in local areas. Since the beginning of 2000, several local dependence measures have been introduced. One of the most important local dependence measures based on the conditional expectation concept was introduced by Bairamov-Kotz (2000) and further developed by Bairamov et al. (2003). This local dependence measure is a bivariate function, having important properties that allow the description of the degree of dependence in particular points and local areas. In many applications, the local dependence function can determine the dependence between important features at specific points. For example, in medical research, it is very important to determine the dependence between the complex of symptoms in healthy tissue and in infected tissue. In this work, we discuss bivariate local dependence functions and introduce a three-variate local dependence function based on the conditional expectation concept. The main idea is to replace the expected value of one random variable with its best predictor expressed in terms of other random variables at some point. Properties of the three-variate local dependence function are given. An example of a three-variate local dependence function with a three-variate normal distribution is also provided. The graphs and tables for the selected parameters are provided. The multivariate extension of the local dependence function is also discussed.         
\end{conclusion}

\bigskip\

\end{document}